
\magnification 1100

  
  %
  %

  %
  %

  \catcode`@=12 

 \def\defrefnote#1{\definexref{#1}{{\the\footnotenumber}}{refnotes}}

  %
  %


 \input eplain.tex
\makeatletter
\def\numberedfootnote{%
ÊÊ\global\advance\footnotenumber by 1
ÊÊ\@eplainfootnote{{\number\footnotenumber}}%
}%
\def\makecolumns#1/#2 {\par \begingroup
ÊÊ \@columndepth = #1
ÊÊ \advance\@columndepth by -1
ÊÊ \divide \@columndepth by #2
ÊÊ \advance\@columndepth by 1
ÊÊ \@linestogoincolumn = \@columndepth
ÊÊ \@linestogo = #1
ÊÊ \currentcolumn = 1
ÊÊ \def\@endcolumnactions{%
ÊÊÊÊÊÊ\ifnum \@linestogo<2
ÊÊÊÊÊÊÊÊ \the\crtok \egroup \endgroup \par 
ÊÊÊÊÊÊ\else
ÊÊÊÊÊÊÊÊ \global\advance\@linestogo by -1
ÊÊÊÊÊÊÊÊ \ifnum\@linestogoincolumn<2
ÊÊÊÊÊÊÊÊÊÊÊÊ\global\advance\currentcolumn by 1
ÊÊÊÊÊÊÊÊÊÊÊÊ\global\@linestogoincolumn = \@columndepth
ÊÊÊÊÊÊÊÊÊÊÊÊ\the\crtok
ÊÊÊÊÊÊÊÊ \else
ÊÊÊÊÊÊÊÊÊÊÊÊ&\global\advance\@linestogoincolumn by -1
ÊÊÊÊÊÊÊÊ \fi
ÊÊÊÊÊÊ\fi
ÊÊ }%
ÊÊ \makeactive\^^M
ÊÊ \letreturn \@endcolumnactions
ÊÊ \@columnwidth = \hsize
ÊÊÊÊ \advance\@columnwidth by -\parindent
ÊÊÊÊ \divide\@columnwidth by #2
ÊÊ \penalty\abovecolumnspenalty
ÊÊ \noindent 
ÊÊ \valign\bgroup
ÊÊÊÊ &\hbox to \@columnwidth{\strut \hsize = \@columnwidth ##\hfil}\cr
}%
\makeatother

\lefteqnumbers
   \def\testd{oui}
   \def\choixlat{\ifx\numadroite\testd\righteqnumbers
            \else  \lefteqnumbers\fi}
    \choixlat

\catcode`@=\letter
\def\@eplainfootnote#1{\let\@sf\empty 
  \ifhmode\edef\@sf{\spacefactor\the\spacefactor}\/\fi
  \global\advance\hlfootlabelnumber by 1
  \hlstart@impl{foot}{\hlfootlabel}%
  \hldest@impl{footback}{\hlfootbacklabel}%
  \hbox{$^{(#1)}$}%
  \hlend@impl{foot}%
  \@sf\vfootnote{#1.}%
}%
\catcode`@=\other

  \interfootnoteskip=0pt
  
  \everyfootnote={\eightpoint\leftskip=5truemm\rightskip5truemm}
  
  \hsize150truemm\vsize 240truemm\hoffset=5truemm
  \def\dimstand{\hsize 150truemm\vsize 240truemm\hoffset=5truemm\voffset=0truemm}
  
  \pretolerance=500\tolerance=1000\brokenpenalty=5000
  \parindent3mm
  
  \countdef\temps=170
  \temps=\time
  \countdef\nminutes=171{\nminutes = \time}
  \countdef\nheure=172
  \def\heure{\begingroup                     
     \temps = \time \divide\temps by 60
     \nheure = \temps                        
     \nminutes = \time
     \multiply\temps by 60
     \advance\nminutes by -\temps            
     \ifnum\nminutes<10 \toks1 = {0}%
     \else\toks1 = {}%
     \fi
     \number\nheure h\the\toks1 \number\nminutes  
  \endgroup}%

  \newcount\chstart
  \chstart=\pageno
 \headline={\ifnum\pageno=\chstart {\hfill} \else {\hss \tenrm --\ \folio\ --\hss}\fi}
  \footline={\hfill}
  \normalbaselines
  \frenchspacing
    \def\dater{\vglue-10mm\rightline{(\the\day/\the\month/\the\year)}}
  \def\dateheure{\vglue-10mm\rightline{(\the\day/\the\month/\the\year,\ \heure)}}

  \newif\ifpagetitre \pagetitretrue
\newtoks\hautpagetitre \hautpagetitre={\hfill}
\newtoks\baspagetitre \baspagetitre={\hfill}
\newtoks\auteurcourant \auteurcourant={\hfill}
\newtoks\titrecourant \titrecourant={\hfill}
\newtoks\hautpagegauche
\newtoks\hautpagedroite
\newtoks\hautpagemilieu
\hautpagemilieu={\tenrm\hfil -- \folio\ -- \hfil}
\hautpagegauche={\ifx\midfolio\oui\the\hautpagemilieu\else\tenrm\folio\hfill\the\auteurcourant\hfill\fi}
\hautpagedroite={\ifx\midfolio\oui\the\hautpagemilieu\else\hfill\the\titrecourant\hfill\tenrm\folio\fi}
\newtoks\baspagegauche \baspagegauche={\hfil}
\newtoks\baspagedroite \baspagedroite={\hfil}
\headline={\ifpagetitre\the\hautpagetitre
\else\ifodd\pageno\the\hautpagedroite\else\the\hautpagegauche\fi\fi }
\footline={\ifpagetitre\the\baspagetitre
\else\ifodd\pageno\the\baspagedroite
\else\the\baspagegauche\fi\fi \global\pagetitrefalse}

\def\pageblanche{\vfill\eject\pagetitretrue
\null\vfill\eject
\pagetitretrue
}
\def\chgtpage{\ifodd\pageno \else
\pageblanche \fi \pagetitretrue\titreun=0\footnotenumber=0}

\def\chgtpageincrtitreun{\ifodd\pageno \else
\pageblanche \fi \pagetitretrue\footnotenumber=0}

\def\majnombres{\ifodd\pageno \else
\pageblanche \fi \pagetitretrue\hautpoly\titreun=0\footnotenumber=0}

\def\hautspages#1#2{\auteurcourant={\ninepcap#1}\titrecourant={\nineit#2}}

\ifnum\chstart=\pageno \pagetitretrue\fi
  

  \def\PAR{\par}
  
  \def\leftnote#1{\vadjust{\setbox1=\vtop{\hsize 20mm\parindent=0pt\eightpoint
  \baselineskip=9pt\rightskip=4mm plus 4mm\vskip-4.7mm#1}\hbox{\kern-2cm\smash{\box1}}}}

  
  \def\raggedcenter{\leftskip=20pt plus 10em  
       \rightskip=\leftskip 
        \parfillskip=0pt 
         \spaceskip=.3333em \xspaceskip=.5em 
          \pretolerance=9999 \tolerance=9999
           \hyphenpenalty=9999 \exhyphenpenalty=9999 }
           
  \def\titrecentre#1{{\parindent0mm\raggedcenter
       \spaceskip=.6em plus .2em minus .2em\xspaceskip=.6em plus .2em minus .2em
        \tit#1\par}}
        


  \def\oui{oui}
  
   \def\fontetitreun{\ifx\paradouze\oui\douzepts\gpdouze\twelvebf\textfont1=\twelveib\else
\quatorzepts\gpquatorze\fourteenbf\fi}

\def\fontetitreunl{\douzepts\textfont1=\twelveib\scriptfont1=\tenib\fourteenti}
 
 \def\fontetitredeux{\textfont1=\eleveni\ifx\paradouze\oui\onzepts\scriptfont1=\ninei\elevenit\else
                        \douzepts\twelveit\fi}
 
   \def\fontetitredeuxb{\ifx\paradouze\oui\onzepts\eleventi\gponze\textfont1=\elevenib\scriptfont1=\nineib
                         \else\douzepts\twelveti\scriptfont1=\twelveib\scriptfont1=\tenib\gpdouze\fi}
                         
\def\fontetitredeuxl{\onzepts\textfont1=\elevenbf\scriptfont1=\ninebf\twelvebf}
  
\def\fontetitretrois{\textfont0=\elevenrm\scriptfont0=\eightrm\textfont1=\eleveni
                      \scriptfont1=\eighti\scriptscriptfont1=\sixi\elevenit}
                      
\def\fontetitrequatre{\textfont0=\elevenrm\scriptfont0=\eightrm\textfont1=\eleveni
                      \scriptfont1=\eighti\scriptscriptfont1=\sixi\elevenrm}
  
  \newcount\titreun\titreun=0
  \newcount\titredeux\titredeux=0
  \newcount\titretrois\titretrois=0
  \newcount\titrequatre\titrequatre=0
  \newcount\enonce\enonce=0
  
  \def\incr#1{\global\advance#1 by 1 {\the #1}}
  \def\avance#1{\global\advance#1 by 1}
  \def\init#1{\global#1=0}
  
  \long\def\Indentation#1#2{\setbox10=\hbox{\fontetitreun#1}
  	                    \ifdim\wd10 < 4mm
                         \setbox10=\hbox to 4mm{\box10\hfill}
                       \else\ifdim\wd10 < 6mm
                         \setbox10=\hbox to 6mm{\box10\hfill}
  	                    \else\ifdim\wd10 < 8mm
                         \setbox10=\hbox to 8mm{\box10\hfill}
                       \else\ifdim\wd10 < 12mm
                         \setbox10=\hbox to 12mm{\box10\hfill}
                       \fi\fi\fi\fi
                       \dimen10=\hsize
                       \advance \dimen10 by -\wd10
                       \noindent \box10 %
                       \ignorespaces
                       \hbox{\vtop{\hsize=\dimen10\raggedright\noindent\fontetitreun#2}}}

  \long\def\paraun#1{\removelastskip\par\medskip\goodbreak\vskip0pt plus.01\vsize\penalty-100
                \vskip0pt plus-.01\vsize
  	              \init{\titredeux}\ifnum\optionparag=1{\init\eqnumber\init\enonce}\else{}\fi
                  \goodbreak{\fontetitreun
  	                \Indentation{\incr{\titreun}.\ }{\fontetitreun #1\par}}\nobreak\medskip}

 %
 %
 \long\def\paraunc#1{\removelastskip\par\bigskip\goodbreak\vskip0pt plus.01\vsize\penalty-100
                \vskip0pt plus-.01\vsize
  	              \init{\titredeux}
                 \ifnum\optionparag=1{\init{\eqnumber}\init\enonce}\else{}\fi
                  \goodbreak
  	                {\parindent0mm\raggedcenter\fontetitreun\incr{\titreun}.\ 
                     \fontetitreun #1\par}\nobreak\medskip}
                     
\newtoks\titreunl
\titreunl={\ifnum\titreun=1{I}\fi%
\ifnum\titreun=2{II}\fi%
\ifnum\titreun=3{III}\fi%
\ifnum\titreun=4{IV}\fi%
\ifnum\titreun=5{V}\fi%
\ifnum\titreun=6{VI}\fi%
\ifnum\titreun=7{VII}\fi%
\ifnum\titreun=8{VIII}\fi%
\ifnum\titreun=9{IX}\fi%
\ifnum\titreun=10{X}\fi%
\ifnum\titreun=11{XI}\fi%
\ifnum\titreun=12{XII}\fi%
\ifnum\titreun=13{XIII}\fi%
}
\long\def\paraunl#1{\removelastskip\par\bigskip\bigskip\goodbreak\vskip0pt plus.01\vsize\penalty-100
                \vskip0pt plus-.01\vsize
  	              \init{\titredeux}\ifnum\optionparag=1{\init\eqnumber\init\enonce}\else{}\fi
                  \goodbreak{\fontetitreunl
  	                \Indentation{\global\advance\titreun by 1{\the\titreunl}.\ }{\fontetitreunl #1\par}}\nobreak\smallskip}

  
  \long\def\paradeux#1{\init{\titretrois}\vskip0pt plus.01\vsize\penalty-10
                \vskip0pt plus-.01\vsize\ifx \elie\oui\medskip\ifnum\titredeux>0\medskip\fi\fi
                 \Indentation{\fontetitredeux\the\titreun${\cdot}$\incr{\titredeux}.}
                              {\fontetitredeux\textfont1=\eleveni#1}\nobreak\par }
  
  \long\def\paradeuxb#1{\init{\titretrois}\vskip0pt plus.001\vsize\penalty-10
                \vskip0pt plus-.01\vsize{\ifx \elie\oui\medskip\ifnum\titredeux>0\medskip\fi\fi
                  \Indentation
  {\fontetitredeuxb\the\titreun${\cdot}$\incr{\titredeux}.}{ \fontetitredeuxb#1}}\nobreak
\smallskip}

\newtoks\titredeuxl
\titredeuxl={\ifnum\titredeux=1{A}\fi%
\ifnum\titredeux=2{B}\fi%
\ifnum\titredeux=3{C}\fi%
\ifnum\titredeux=4{D}\fi%
\ifnum\titredeux=5{E}\fi%
\ifnum\titredeux=6{F}\fi%
\ifnum\titredeux=7{G}\fi%
\ifnum\titredeux=8{H}\fi%
\ifnum\titredeux=9{I}\fi%
\ifnum\titredeux=10{J}\fi%
\ifnum\titredeux=11{K}\fi%
\ifnum\titredeux=12{L}\fi%
\ifnum\titredeux=13{M}\fi%
}
 \long\def\paradeuxl#1{\init{\titretrois}\vskip0pt plus.001\vsize\penalty-10
                \vskip0pt plus-.01
                \vsize \bigskip%
                  \Indentation
     {\fontetitredeuxl\global\advance\titredeux by 1
  \quad \the\titreunl${\cdot}$\the\titredeuxl.}{ \fontetitredeuxl#1}
  \removelastskip\nobreak\smallskip}
  

  \long\def\paratrois#1{\init{\titrequatre}\ifdim\lastskip<\smallskipamount
                \removelastskip\smallskip\fi
                 \vskip0pt plus.01\vsize\penalty-10
                  \vskip0pt
plus-.01\vsize{\ifx \elie\oui\ifnum\titretrois>0\medskip\fi\fi
\Indentation{\fontetitretrois\the\titreun${\cdot}$\the\titredeux${\cdot}$\incr{\titretrois}.\ }
  {\hskip0mm\baselineskip=14pt\fontetitretrois#1}\nobreak\smallskip}}
  
  
  \long\def\paratroisl#1{\init{\titrequatre}\ifdim\lastskip<\smallskipamount
                \removelastskip\fi
                 \vskip0pt plus.01\vsize\penalty-10
                  \vskip0pt
plus-.01\vsize\ifx \elie\oui\bigskip
\fi
\Indentation{\fontetitretrois\quad \quad \the\titreunl{${\cdot}$}\the\titredeuxl${\cdot}$\incr{\titretrois}.\ }
  {\hskip0mm\fontetitretrois#1}\nobreak\smallskip}


  \long\def\paraquatre#1{\ifdim\lastskip<\smallskipamount
                \removelastskip\smallskip\fi
                 \vskip0pt plus.01\vsize\penalty-10
                  \vskip0pt
                  plus-.01\vsize\par
 
\Indentation{\fontetitrequatre \the\titreun{${\cdot}$}\the\titredeux${\cdot}$\the\titretrois${\cdot}$\incr{\titrequatre}.\ }
{\hskip0mm\fontetitrequatre#1}\nobreak\smallskip}


\newtoks\titrequatrel
\titrequatrel={\ifnum\titrequatre=1{a}\fi%
\ifnum\titrequatre=2{b}\fi%
\ifnum\titrequatre=3{c}\fi%
\ifnum\titrequatre=4{d}\fi%
\ifnum\titrequatre=5{e}\fi%
\ifnum\titrequatre=6{f}\fi%
\ifnum\titrequatre=7{g}\fi%
\ifnum\titrequatre=8{h}\fi%
\ifnum\titrequatre=9{i}\fi%
\ifnum\titrequatre=10{j}\fi%
\ifnum\titrequatre=11{k}\fi%
\ifnum\titrequatre=12{l}\fi%
\ifnum\titrequatre=13{m}\fi%
}
\long\def\paraquatrel#1{\ifdim\lastskip<\smallskipamount
                \removelastskip\smallskip\fi
                 \vskip0pt plus.01\vsize\penalty-10
                  \vskip0pt
                  plus-.01\vsize{\bigskip
\Indentation{\global\advance\titrequatre by 1
\fontetitrequatre\quad \quad \quad \the\titreunl${\cdot}$\the\titredeuxl${\cdot}$\the\titretrois${\cdot}$\the\titrequatrel.\ }
{\hskip0mm\fontetitrequatre#1}\nobreak\smallskip}}

\ifx\optionkeys\oui
\def\drefun#1{\definexref{¤#1}{{\the\titreun}}{}} 
\def\drefdeux#1{\definexref{¤#1}{{\the\titreun}.{\the\titredeux}}{}}
\def\dreftrois#1{\definexref{¤#1}{{\the\titreun}.{\the\titredeux}.{\the\titretrois}}{}}
\else
\def\drefun#1{\definexref{prg#1}{{\the\titreun}}{}} 
\def\drefdeux#1{\definexref{prg#1}{{\the\titreun}.{\the\titredeux}}{}}
\def\dreftrois#1{\definexref{prg#1}{{\the\titreun}.{\the\titredeux}.{\the\titretrois}}{}}
\fi

%


  \long\def\propdeux#1#2#3#4{%
       \avance{\enonce}
       \leavevmode\edef\temp{#2}%
         \ifx\temp\empty 
          \else
           \definexref{#2}{#1~{\the\titreun.\the\enonce}}{enonces}
            \definexref{s#2}{{\the\titreun.\the\enonce}}{enonces}
             \fi
\smallskip
      \noindent{\bf#1\ {\bf\the\titreun.\the\enonce{#3}.}\enspace}{\sl#4\par}%
      \ifdim\lastskip<\medskipamount \removelastskip\penalty55\par \fi
   }

  \long\def\propun#1#2#3#4{%
      \avance{\enonce}
       \leavevmode\edef\temp{#2}%
        \ifx\temp\empty 
          \else
           \definexref{#2}{#1~{\the\enonce}}{enonces}
            \definexref{{s#2}}{{\the\enonce}}{enonces}
             \fi
   \par 
     \noindent{\bf#1\ {\bf\the\enonce{#3}.}\enspace}{\sl#4\par}%
     \ifdim\lastskip<\medskipamount \removelastskip\penalty55\medskip\fi
  }
  
  \long\def\prop#1#2#3#4{\ifnum\optionparag=1
                          \propdeux{#1}{#2}{\textfont1=\elevenib#3}{#4} \else\propun{#1}{#2}{\textfont1=\elevenib#3}{#4}\fi}

  \long\def\propt#1#2#3{\ifx\tpf\oui \prop{Th\'eo\-r\`eme}{#1}{#2}{#3}\par
                       \else\prop{Theorem}{#1}{#2}{#3}\par\fi}
  \long\def\Propt#1#2{\propt{#1}{}{#2}}
  \long\def\propl#1#2#3{\ifx\tpf\oui\prop{Lem\-me}{#1}{#2}{#3}\par
                         \else\prop{Lemma}{#1}{#2}{#3}\par\fi}
  \long\def\Propl#1#2{\propl{#1}{}{#2}}
  \long\def\propc#1#2#3{\ifx\tpf\oui\prop{Corol\-laire}{#1}{#2}{#3}\par
                         \else\prop{Corollary}{#1}{#2}{#3}\par\fi}

  \long\def\propd#1#2#3{\ifx\tpf\oui\prop{D\'efi\-nition}{#1}{#2}{#3}\par
                       \else\prop{Definition}{#1}{#2}{#3}\par\fi} 
  
  \long\def\proptd#1#2#3{\ifx\tpf\oui\prop{Th\'eor\`eme et d\'efi\-nition}{#1}{#2}{#3}\par
                       \else\prop{Theorem and definition}{#1}{#2}{#3}\par\fi}


  
  \newcount\optionparag\optionparag=1
  
  \long\def\section#1#2{\ifnum\optionparag=1 \paraun{#2} 
                        \else\goodbreak{\fontetitreun
  	                \Indentation{#1.\ }{#2}}\nobreak\smallskip\fi}

  \def\eqconstruct#1{\ifnum\optionparag=1{\the\titreun\hbox{$\cdot$}#1}\else{#1}\fi}

  
  
  \def\numref{oui}  
  
  \newcount\mesref\mesref=0 
  \def\defbib#1{\ifx\numref\oui\global\advance\mesref by 1 \definexref{#1}{{\the
                 \mesref}}{}\else\definexref{#1}{#1}{}\fi}
  \def\bibtem#1{\defbib{#1}\item{\citer{#1}}}
  \def\citer#1{[\ref{#1}]}
  \def\citeplus#1#2{[\ref{#1}; #2]}

  
  \font\seventeenmsa=msam10 at 17pt    
  \font\fourteenmsa=msam10 at 14pt
  \font\twelvemsa=msam10 at 12pt
  \font\tenmsa=msam10                 
  \font\ninemsa=msam10 at 9pt 
  \font\eightmsa=msam10 at 8pt 
  \font\sevenmsa=msam7 
  \font\sixmsa=msam10 at 6pt
  \font\fivemsa=msam5
  \newfam\msafam\textfont\msafam=\tenmsa\scriptfont\msafam=\sevenmsa\scriptscriptfont\msafam=\fivemsa
  
  \font\seventeenbb=msbm10 at 17pt     
  \font\fourteenbb=msbm10 at 14pt
  \font\twelvebb=msbm10 at 12pt
  \font\tenbb=msbm10                   
  \font\ninebb=msbm10 at 9pt 
  \font\eightbb=msbm10 at 8pt 
  \font\sevenbb=msbm7 
  \font\sixbb=msbm10 at 6pt
  \font\fivebb=msbm5 
  \newfam\bbfam\textfont\bbfam=\tenbb\scriptfont\bbfam=\sevenbb\scriptscriptfont\bbfam=\fivebb
  \def\bb{\fam\bbfam\tenbb}%

  \font\seventeenscaln=eusm10 at 17pt   
  \font\twelvescaln=eusm10 at 12pt
  \font\tenscaln=eusm10                
  \font\ninescaln=eusm10 scaled 900
  \font\eightscaln=eusm10 scaled 800
  \font\sevenscaln=eusm10 scaled 700
  \font\sixscaln=eusm10 scaled 600
   
  \newfam\scalnfam\textfont\scalnfam=\tenscaln\scriptfont\scalnfam=\sevenscaln\scriptscriptfont\scalnfam=\sixscaln
  \def\scaln{\fam\scalnfam\tenscaln}%
  \def\scal{\scaln}
  
  \font\tenscalb=eusb10                

  \font\sevenscalb=eusb10 scaled 700

  \newfam\scalbfam\textfont\scalbfam=\tenscalb\scriptfont\scalbfam=\sevenscalb
  %
  
  %
  %
  \font\fourteenrm=cmr12 scaled 1200
  \font\elevenrm=cmr10 at 11pt
  \font\twelverm=cmr12
  \font\ninerm=cmr9
  \font\eightrm=cmr8      
  \font\sevenrm=cmr7
  \font\sixrm=cmr6

  \font\seventeenpcap=cmcsc10 at 17pt
  \font\tenpcap=cmcsc10                        
  \font\ninepcap=cmcsc9
  \font\eightpcap=cmcsc8
  \font\sevenpcap=cmcsc10 scaled 700
  
  \newfam\pcapfam\textfont\pcapfam=\tenpcap\scriptfont\pcapfam=\sevenpcap
  \def\pcap{\fam\pcapfam\tenpcap}
  
  \font\seventeenrm=cmbx12 scaled 1400

  \font\fourteenbf=cmbx10 scaled 1400
  
  \font\twelvebf=cmbx12
  \font\elevenbf=cmbx10 at 11pt
  \font\ninebf=cmbx9  
  \font\eightbf=cmbx8
  \font\sixbf=cmbx6
  
  \font\tengot=eufm10                           
   
  \font\eightgot=eufm10 at 8truept 
  \font\sevengot=eufm7 
  \font\sixgot=eufm10 at 6 truept 
   
  \newfam\gotfam
  \textfont\gotfam=\tengot\scriptfont\gotfam=\sevengot\scriptscriptfont\gotfam=\sixgot
  \def\got{\fam\gotfam\tengot}%

  
  \def\tit{%
  \textfont0=\seventeenrm\scriptfont0=\tenrm\def\rm{\fam0\seventeenrm}%
  \textfont1=\seventeenib\scriptfont1=\twelveib%
  \textfont2=\seventeensy\scriptfont2=\twelvesy\scriptscriptfont2=\ninesy
  \textfont3=\seventeenex
  \textfont\itfam=\seventeenti
  \def\it{\fam\itfam\seventeenti}%
  \textfont\bbfam=\seventeenbb \scriptfont\bbfam=\twelvebb
  \def\bb{\fam\bbfam\seventeenbb}%
  \textfont\msafam=\seventeenmsa\scriptfont\msafam=\twelvemsa
  \textfont\scalnfam=\seventeenscaln
  \def\pcap{\fam\pcapfam\seventeenpcap}
  \normalbaselineskip=25pt\normalbaselines\rm}

  \font\seventeenti=cmbxti10 scaled 1680
  
  \font\fourteenti=cmbxti10 at 14pt
  
  \font\twelveti=cmbxti10 scaled 1200
  \font\eleventi=cmbxti10 at 11pt

  %
  %
  \font\twelveit=cmti12	
  \font\elevenit=cmti10 scaled 1100
  \font\nineit=cmti9
  \font\eightit=cmti8
  \font\sevenit=cmti7

  %
  %
 
 \font\seventeenib=cmmib10 scaled 1680
  \font\fourteenib=cmmib10 scaled 1400
  \font\twelveib=cmmib10 scaled 1200
  \font\elevenib=cmmib10 scaled 1100
  \font\tenib=cmmib10
\font\eightib=cmmib10 scaled 800
  \font\nineib=cmmib10 scaled 900
\font\sevenib=cmmib10 scaled 700
\font\sixib=cmmib10 scaled 600
\font\fiveib=cmmib10 scaled 500

\ifx\ITAN\oui
\else
\innernewfam\cmmibfam
\textfont\cmmibfam=\tenib
\scriptfont\cmmibfam=\sevenib
\scriptscriptfont\cmmibfam=\fiveib
\def\ib{\fam\cmmibfam\tenib}
\fi

  %
  %
  
  \font\eleveni=cmmi10 scaled 1100
  \font\ninei=cmmi9
  \font\eighti=cmmi8 
  \font\seveni=cmmi7 			                
  \font\sixi=cmmi6
  
  \font\ninesl=cmsl9                    
  \font\eightsl=cmsl8 
  \font\sevensl=cmsl10 at 7pt

  \font\ninett=cmtt9                    
  \font\eighttt=cmtt8
  \font\seventt=cmtt10 scaled 700

  \font\seventeensy=cmsy10 scaled 1680    
  \font\fourteensy=cmsy10 scaled 1400
  \font\twelvesy=cmsy10 scaled 1176
  
  \font\ninesy=cmsy9                      
  \font\eightsy=cmsy8
  \font\sixsy=cmsy6
  \font\seventeenex=cmex10 at 17pt
  \font\fourteenex=cmex10 at 14pt
  \font\twelveex=cmex10 at 12pt
  \font\nineex=cmex10 at 9pt
  \font\eightex=cmex10 at 8pt
  \font\sevenex=cmex10 at 7pt
  \font\sixex=cmex10 at 6pt
  \font\fiveex=cmex10 at 5pt
  
   
  \font\fourteengp=cmmi10 at 14pt
  
  \font\twelvegp=cmmib10 at 12pt
  \font\elevengp=cmmib10 at 11pt
  \font\tengp=cmmib10                          
  \font\ninegp=cmmib10 at 9pt 
  \font\eightgp=cmmib8 
   
  \font\sixgp=cmmib6


  \def\gponze{\textfont0=\elevenbf\scriptfont0=\eightbf\scriptscriptfont0=\sixbf
           \textfont1=\elevengp\scriptfont1=\eightgp\scriptscriptfont1=\sixgp}
  \def\gpdouze{\textfont0=\twelvebf\scriptfont0=\tenbf\scriptscriptfont0=\ninebf
           \textfont1=\twelvegp\scriptfont1=\tengp\scriptscriptfont1=\ninegp}        
  
 \def\gpquatorze{\textfont0=\fourteenbf\scriptfont0=\twelvebf\scriptscriptfont0=\elevenbf
           \textfont1=\fourteengp\scriptfont1=\twelvegp\scriptscriptfont1=\elevengp}

  
  \expandafter\chardef\csname pre amssym.def at\endcsname=\the\catcode`\@
  \catcode`\@=11
  \def\undefine#1{\let#1\undefined}
  \def\newsymbol#1#2#3#4#5{\let\next@\relax
   \ifnum#2=\@ne\let\next@\msafam@\else
   \ifnum#2=\tw@\let\next@\bbfam@\fi\fi
   \mathchardef#1="#3\next@#4#5}
  \def\mathhexbox@#1#2#3{\relax
   \ifmmode\mathpalette{}{\m@th\mathchar"#1#2#3}%
   \else\leavevmode\hbox{$\m@th\mathchar"#1#2#3$}\fi}
  \def\hexnumber@#1{\ifcase#1 0\or 1\or 2\or 3\or 4\or 5\or 6\or 7\or 8\or
   9\or A\or B\or C\or D\or E\or F\fi}
  
  \def\setboxz@h{\setbox\z@\hbox}
  \def\wdz@{\wd\z@}
  \def\boxz@{\box\z@}
  
  \edef\msafam@{\hexnumber@\msafam}
  \mathchardef\dabar@"0\msafam@39
  
  \edef\bbfam@{\hexnumber@\bbfam}
  \def\widehat#1{\setboxz@h{$\m@th#1$}%
   \ifdim\wdz@>\tw@ em\mathaccent"0\bbfam@5B{#1}%
   \else\mathaccent"0362{#1}\fi}
  \def\widetilde#1{\setboxz@h{$\m@th#1$}%
   \ifdim\wdz@>\tw@ em\mathaccent"0\bbfam@5D{#1}%
   \else\mathaccent"0365{#1}\fi}
  \newsymbol\leqq 1335          
  \newsymbol\leqslant 1336
  \newsymbol\lessgtr 1337       
  \newsymbol\backprime 1038     
  \newsymbol\risingdotseq 133A  
  \newsymbol\fallingdotseq 133B 
  \newsymbol\succcurlyeq 133C   
  \newsymbol\geqq 133D          
  \newsymbol\geqslant 133E
  \newsymbol\nmid 232D
  \newsymbol\nexists 2040
  \newsymbol\smallsetminus 2272
  \newsymbol\varnothing 203F
  
  \catcode`\@=\active

  \catcode`\@=11
  \newcount\typofr\typofr=1
  
  \catcode`\;=\active
  \def;{\ifnum\typofr=1\relax\ifhmode\ifdim\lastskip>\z@\unskip\fi
     \kern.2em\fi\string;\else\string;\fi}
  
  \catcode`\:=\active
  \def:{\ifnum\typofr=1\relax\ifhmode\ifdim\lastskip>\z@\unskip\fi
  \penalty\@M\ \fi\string:\else\string:\fi}
  
  \catcode`\!=\active
  \def!{\ifnum\typofr=1\relax\ifhmode\ifdim\lastskip>\z@\unskip\fi
     \kern.2em\fi\string!\else\string!\fi}
  
  \catcode`\?=\active
  \def?{\ifnum\typofr=1\relax\ifhmode\ifdim\lastskip>\z@\unskip\fi
     \kern.2em\fi\string?\else\string?\fi}

  \def\francais{\typofr=1\def\tpf{oui}}
  \def\anglais{\typofr=2\def\tpf{non}\def\english{oui}}
  \def\oui{oui}
  \francais
  
  \catcode`\@=12
  

\ifx\textures\oui
\def\raye #1|{\leavevmode\setbox1=\hbox{#1}%
\raise .5pt\hbox to \wd1{\xleaders\hbox{\rge{$ \sct / $}%
\kern 1pt}\hfill\kern -1pt }\kern -\wd1 \unhbox1\relax }

\def\barre #1|{\leavevmode\setbox1=\hbox{#1}%
\rlap{\Red\vrule height 2.4pt depth -1.2pt width \wd1}\Black \unhbox1\relax}
\else
\def\raye #1|{\leavevmode\setbox1=\hbox{#1}%
\raise .5pt\hbox to \wd1{\xleaders\hbox{\rge{$ \sct / $}%
\kern 1pt}\hfill\kern -1pt }\kern -\wd1 \unhbox1\relax }

\def\barre #1|{\leavevmode\setbox1=\hbox{#1}%
\rlap{\color{red}\vrule height 2.4pt depth -1.2pt width \wd1}\color{black} \unhbox1\relax}

\fi
  

  
  \def\og{\leavevmode\raise.24ex\hbox{$\scriptscriptstyle\langle\!\langle\>$}}    
  \def\fg{\leavevmode\raise.24ex\hbox{$\scriptscriptstyle\>\rangle\!\rangle$}}    

  \def\d{\,{\rm d}}
  \def\dd{{\rm d}}

  \def\r{{\bb R}}
  
  \def\N{{\bb N}}

  \def\E{{\scal E}}

  \def\M{{\scal M}}
  
  \def\O{{\scal O}}
  \def\P{{\scaln P}}

  \def\frac#1#2{{#1\over #2}}
  \def\di#1#2{\sct#1\atop{\sct#2}}

  \def\numero{n$^{\rm o}\thinspace$}
\def\numeros{n$^{\rm os}\thinspace$}

  \def\qedbox{$\rlap{$\sqcap$}\sqcup$}           
  \def\qed{\nobreak\hfill\penalty250 \hbox{}\nobreak\hfill\qedbox\par }

  \def\numero{n$^{\rm o}\thinspace$}

  \def\pnt{prime number theorem}

  \def\¤{\S\thinspace}

  \def\¥{$\bullet$ }
  
  
  \def\e{{\rm e}}

  \def\no#1{\Vert#1\Vert}
  
  \def\NO#1{\left\Vert#1\right\Vert}

  \def\epsilon{\varepsilon}

  \def\phi{\varphi}
  \def\theta{\vartheta}
  \def\rho{\varrho}
  \def\dm{{\textstyle{1\over 2}}}
  \def\txt{\textstyle}
  \def\dsp{\displaystyle}
  \def\sct{\scriptstyle}
  \def\pf{\noi{\it Proof. }}
  \def\nid{\ifnum\typofr=1\par\noindent{\it D\'emonstration. }\else\pf\fi}
  \def\noi{\noindent}
  \def\rem{\ifnum\typofr=1\noi{\it Remarque.}\ \else\noi{\it Remark.}\ \fi}
  \def\rems{\ifnum\typofr=1\noi{\it Remarques.}\ \else\noi{\it Remarks.}\ \fi}
  \def\re{{\Re e\,}}

  \def\1{{\bf 1}}
  \def\|{\Vert}

  \def\leq{\leqslant}
  \def\geq{\geqslant}

  \def\ie{{i.e.\ }}
  \def\eg{{e.g.}}
  

  \def\fl#1{\left\lfloor #1 \right\rfloor}

  \def\log{\mathop{\rm log}\nolimits}
  \def\ft#1#2{{\txt{#1\over #2}}}


  
\def\vbs#1{\left|#1\right|}
\def\Vbs#1{\bigg|#1\bigg|}
\def\abs#1{\left|#1\right|}


  \def\pmb#1{\setbox0=\hbox{#1}%
  \kern-.025em\copy0\kern-\wd0\kern.05em\copy0\kern-\wd0\kern-.025em\raise .0433em\box0 }

  
  \skewchar\eighti='177 \skewchar\sixi='177
  \skewchar\eightsy='60 \skewchar\sixsy='60
  
  \def\eightpoint{%
  \textfont0=\eightrm\scriptfont0=\sixrm\scriptscriptfont0=\fiverm
  \def\rm{\fam0\eightrm}%
  \textfont1=\eighti\scriptfont1=\sixi
  \scriptscriptfont1=\fivei\def\oldstyle{\fam1\seveni}%
  \textfont2=\eightsy\scriptfont2=\sixsy\scriptscriptfont2=\fivesy
  \textfont3=\eightex\scriptfont3=\sixex
  \textfont\itfam=\eightit
  \def\it{\fam\itfam\eightit}%
  \textfont\slfam=\eightsl
  \def\sl{\fam\slfam\eightsl}%
  \textfont\bbfam=\eightbb \scriptfont\bbfam=\sixbb\scriptscriptfont\bbfam=\fivebb
  \def\bb{\fam\bbfam\eightbb}%
  \textfont\msafam=\eightmsa\scriptfont\msafam=\sixmsa
  \textfont\scalnfam=\eightscaln
  \def\scaln{\fam\scalnfam\eightscaln}
  \textfont\ttfam=\eighttt
  \def\tt{\fam\ttfam\eighttt}%
\textfont\gotfam=\eightgot
  \textfont\bffam=\eightbf\scriptfont\bffam=\sixbf\scriptscriptfont\bffam=\fivebf
  \def\bf{\fam\bffam\eightbf}%
  \ifx\ITAN\oui\else\textfont\cmmibfam=\eightib
       \scriptfont\cmmibfam=\sixib
        \scriptscriptfont\cmmibfam=\fiveib
         \def\ib{\fam\cmmibfam\eightib}
   \fi
  \textfont\pcapfam=\eightpcap
  \def\pcap{\fam\pcapfam\eightpcap}
  \abovedisplayskip=2pt plus2pt minus 2pt
  \belowdisplayskip=2pt plus1pt minus 2pt
  \abovedisplayshortskip= 1pt plus 2pt minus 1pt
  \belowdisplayshortskip= 1pt plus 2pt minus 1pt
  \smallskipamount=2pt plus 1pt minus 2pt
  \medskipamount=3pt plus 2pt minus 2pt
  \bigskipamount=7pt plus 3pt minus 3pt
  \setbox\strutbox=\hbox{\vrule height 5pt depth 2pt width 0pt}%
  \normalbaselineskip=9pt\normalbaselines\rm}

  \def\({\left(}
  \def\){\right)}
  
  \def\footnoterule{\kern -2pt\hrule width 7truecm\kern 2.4pt}
  
  \def\xnotedef#1{\definexref{#1}{\noexpand\number\footnotenumber}{Note}}%

  
  
  \def\ninepoint{%
  \textfont0=\ninerm\scriptfont0=\sixrm\scriptscriptfont0=\fiverm
  \def\rm{\fam0\ninerm}%
  \textfont1=\ninei\scriptfont1=\sixi
  \scriptscriptfont1=\fivei\def\oldstyle{\fam1\ninei}%
  \textfont2=\ninesy\scriptfont2=\sixsy\scriptscriptfont2=\fivesy
  \textfont3=\nineex\scriptfont3=\sixex
  \textfont\itfam=\nineit
  \def\it{\fam\itfam\nineit}%
  \textfont\slfam=\ninesl
  \def\sl{\fam\slfam\ninesl}%
  \textfont\bbfam=\ninebb\scriptfont\bbfam=\sixbb\scriptscriptfont\bbfam=\fivebb
  \def\bb{\fam\bbfam\ninebb}%
  \textfont\msafam=\ninemsa\scriptfont\msafam=\sixmsa\scriptscriptfont\msafam=\fivemsa
  \textfont\scalnfam=\ninescaln
  \def\scaln{\fam\scalnfam\ninescaln}
  \textfont\ttfam=\ninett
  \def\tt{\fam\ttfam\ninett}%
  \textfont\bffam=\ninebf\scriptfont\bffam=\sixbf\scriptscriptfont\bffam=\fivebf
  \def\bf{\fam\bffam\ninebf}%
  \abovedisplayskip=3pt plus2pt minus 2pt
  \belowdisplayskip=3pt plus1pt minus 2pt
  \abovedisplayshortskip= 2pt plus 2pt minus 1pt
  \belowdisplayshortskip= 2pt plus 2pt minus 1pt
  \smallskipamount=2pt plus 1pt minus 2pt
  \medskipamount=3pt plus 2pt minus 2pt
  \bigskipamount=7pt plus 3pt minus 3pt
  \setbox\strutbox=\hbox{\vrule height 5pt depth 2pt width 0pt}%
  \normalbaselineskip=10.5pt plus.3pt minus.3pt\normalbaselines\rm}

  \def\sevenpoint{%
  \textfont0=\sevenrm\scriptfont0=\sixrm\scriptscriptfont0=\fiverm
  \def\rm{\fam0\sevenrm}%
  \textfont1=\seveni\scriptfont1=\sixi
  \scriptscriptfont1=\fivei\def\oldstyle{\fam1\seveni}%
  \textfont2=\sevensy\scriptfont2=\sixsy\scriptscriptfont2=\fivesy
  \textfont3=\sevenex\scriptfont3=\fiveex
  \textfont\itfam=\sevenit
  \def\it{\fam\itfam\sevenit}%
  \textfont\slfam=\sevensl
  \def\sl{\fam\slfam\sevensl}%
  \textfont\bbfam=\sevenbb \scriptfont\bbfam=\sixbb\scriptscriptfont\bbfam=\fivebb
  \def\bb{\fam\bbfam\sevenbb}%
  \textfont\msafam=\sevenmsa\scriptfont\msafam=\sixmsa
  \textfont\scalnfam=\sevenscaln
  \def\scaln{\fam\scalnfam\sevenscaln}
  \textfont\bffam=\sevenbf\scriptfont\bffam=\sixbf\scriptscriptfont\bffam=\fivebf
  \def\bf{\fam\bffam\sevenbf}%
  \textfont\ttfam=\seventt
  \abovedisplayskip=2pt plus2pt minus 2pt
  \belowdisplayskip=2pt plus1pt minus 2pt
  \abovedisplayshortskip= 1pt plus 2pt minus 1pt
  \belowdisplayshortskip= 1pt plus 2pt minus 1pt
  \smallskipamount=2pt plus 1pt minus 2pt
  \medskipamount=3pt plus 2pt minus 2pt
  \bigskipamount=7pt plus 3pt minus 3pt
  \setbox\strutbox=\hbox{\vrule height 5pt depth 2pt width 0pt}%
  \normalbaselineskip=9pt\normalbaselines\rm}

 \def\onzepts{%
 \textfont0=\elevenrm\scriptfont0=\ninerm
 \textfont1=\elevenib\scriptfont1=\ninei}

\def\douzepts{%
  \textfont0=\twelverm\scriptfont0=\tenrm\def\rm{\fam0\twelverm}%
  \textfont1=\twelveib\scriptfont1=\teni%
  \textfont2=\twelvesy\scriptfont2=\tensy\scriptscriptfont2=\eightsy
  \textfont3=\twelveex
  \textfont\itfam=\twelveti
  \def\it{\fam\itfam\twelveti}%
  \textfont\bffam=\twelvebf\scriptfont\bffam=\tenbf\scriptscriptfont\bffam=\eightbf
  \def\bf{\fam\bffam\twelvebf}%
  \textfont\bbfam=\twelvebb \scriptfont\bbfam=\tenbb
  \def\bb{\fam\bbfam\twelvebb}%
  \textfont\msafam=\twelvemsa\scriptfont\msafam=\tenmsa
  \textfont\scalnfam=\twelvescaln
  \normalbaselineskip=15pt\normalbaselines\rm}

\def\quatorzepts{%
  \textfont0=\fourteenrm\scriptfont0=\twelverm\def\rm{\fam0\fourteenrm}%
  \textfont1=\fourteenib\scriptfont1=\twelveib%
  \textfont2=\fourteensy\scriptfont2=\twelvesy\scriptscriptfont2=\tensy
  \textfont3=\fourteenex
  \textfont\itfam=\fourteenti
  \def\it{\fam\itfam\fourteenti}%
  \textfont\bffam=\fourteenbf\scriptfont\bffam=\twelvebf\scriptscriptfont\bffam=\tenbf
  \def\bf{\fam\bffam\fourteenbf}%
  \textfont\bbfam=\fourteenbb \scriptfont\bbfam=\twelvebb
  \def\bb{\fam\bbfam\fourteenbb}%
  \textfont\msafam=\fourteenmsa\scriptfont\msafam=\twelvemsa
  \textfont\scalnfam=\twelvescaln
  \normalbaselineskip=18pt\normalbaselines\rm}


\def\AA{{\it Acta Arith.}}

\def\picture #1 by #2 (#3){\leavevmode\vbox to #2{
     \hrule width #1 height 0pt depth 0pt
      \vfill
       \special{picture #3}}}

\def\illustration #1 by #2 (#3) scaled #4{\dimen1=#2
  \divide\dimen1 by 1000
  \multiply\dimen1 by #4
  \vtop to \dimen1{\dimen1=#1
  \divide\dimen1 by 1000
  \multiply\dimen1 by #4
  \hsize=\dimen1\vss
  \noindent\special{illustration #3 scaled #4}}}

\ifx\optionkeymacros\undefined\else \fi

\catcode`\Œ=\active\defŒ{{\aa}}       
\catcode`\º=\active\defº{\int}        
\catcode`\=\active\def{\c c}        
\catcode`\¶=\active\def¶{\partial}    
\catcode`\Ä=\active\defÄ{\oint}       
\catcode`\Æ=\active\defÆ{\triangle}   
\catcode`\Â=\active\defÂ{\neg}        
\catcode`\µ=\active\defµ{\mu}         
\catcode`\¿=\active\def¿{{\o}}        
\catcode`\¹=\active\def¹{\pi}         
\catcode`\Ï=\active\defÏ{{\oe}}       
\catcode`\§=\active\def§{{\ss}}       
\catcode`\ =\active\def {\dagger}     
\catcode`\Ã=\active\defÃ{\sqrt}       
\catcode`\·=\active\def·{\Sigma}      
\catcode`\Å=\active\defÅ{\approx}     
\catcode`\½=\active\def½{\Omega}      
\catcode`\£=\active\def£{{\it\$}}     
\catcode`\°=\active\def°{\infty}      
\catcode`\¤=\active\def¤{{\S}}        
\catcode`\¦=\active\def¦{{\P}}        
\catcode`\¥=\active\def¥{\bullet}     
\catcode`\»=\active\def»{\leavevmode\raise.585ex\hbox{\b a}}      
\catcode`\¼=\active\def¼{\leavevmode\raise.6ex\hbox{\b o}}        
\catcode`\­=\active\def­{\not=}       
\catcode`\²=\active\def²{\leq}        
\catcode`\³=\active\def³{\geq}        
\catcode`\Ö=\active\defÖ{\div}        
\catcode`\É=\active\defÉ{{\dots}}     
\catcode`\¾=\active\def¾{{\ae}}       
\catcode`\Ç=\active\defÇ{\og}         
\catcode`\Ò=\active\defÒ{``}          
\catcode`\Á=\active\defÁ{!`}          
\catcode`\¢=\active\def¢{\rlap/c}     
\catcode`\Ô=\active\defÔ{`}           
\catcode`\Õ=\active\defÕ{'}           


\catcode`\=\active\def{{\AA}}       
\catcode`\'=\active\def'{\c C}        
\catcode`\¯=\active\def¯{{\O}}        
\catcode`\¸=\active\def¸{\Pi}         
\catcode`\Î=\active\defÎ{{\OE}}       
\catcode`\®=\active\def®{{\AE}}       
\catcode`\×=\active\def×{\diamond}    
\catcode`\¡=\active\def¡{\accent'27}  
\catcode`\Ó=\active\defÓ{''}          
\catcode`\±=\active\def±{\pm}         
\catcode`\È=\active\defÈ{\fg}         
\catcode`\À=\active\defÀ{?`}          
\catcode`\Ð=\active\defÐ{--}          
\catcode`\Ñ=\active\defÑ{---}         


\catcode`\Š=\active\defŠ{\"a}        
\catcode`\'=\active\def'{\"e}        
\catcode`\•=\active\def•{\"{\i}}     
\catcode`\š=\active\defš{\"o}        
\catcode`\Ÿ=\active\defŸ{\"u}        
\catcode`\Ø=\active\defØ{\"y}        
\catcode`\å=\active\defå{\^A}        
\catcode`\€=\active\def€{\"A}        
\catcode`\…=\active\def…{\"O}        
\catcode`\†=\active\def†{\"U}        
\catcode`\‡=\active\def‡{\'a}        
\catcode`\Ž=\active\defŽ{\'e}        
\catcode`\'=\active\def'{\'{\i}}     
\catcode`\—=\active\def—{\'o}        
\catcode`\œ=\active\defœ{\'u}        
\catcode`\ƒ=\active\defƒ{\'E}        
\catcode`\æ=\active\defæ{\^E}        
\catcode`\é=\active\defé{\`E}        %
\catcode`\ˆ=\active\defˆ{\`a}        
\catcode`\=\active\def{\`e}        
\catcode`\"=\active\def"{\`{\i}}     
\catcode`\˜=\active\def˜{\`o}        
\catcode`\=\active\def{\`u}        
\catcode`\Ë=\active\defË{\`A}        
\catcode`\‹=\active\def‹{\~a}        
\catcode`\–=\active\def–{\~n}        
\catcode`\›=\active\def›{\~o}        
\catcode`\Ì=\active\defÌ{\~A}        
\catcode`\"=\active\def"{\~N}        
\catcode`\Í=\active\defÍ{\~O}        
\catcode`\‰=\active\def‰{\^a}        
\catcode`\=\active\def{\^e}        
\catcode`\"=\active\def"{\^{\i}}     
\catcode`\™=\active\def™{\^o}        
\catcode`\ž=\active\defž{\^u}        

\let\optionkeymacros\null

\ifx\montrerlabels\oui
\input montrerlabels.tex
\fi
\newcount\numchapi\numchapi=0

\optionparag=1
\def\paradouze{oui}
\anglais
\def\pnu{p^\nu}
\def\fr{distribution function}
\def\ga{{\got a}}
\def\gb{{\got b}}
\def\gc{{\got c}}
\def\gh{{\got h}}
\def\gR{{\got R}}
\def\FF{{\scal F}}

\dimstand\hautspages{G. Tenenbaum \& J. Verwee}{Effective Erd\H os--Wintner theorems}
\dateheure

\titrecentre{Effective Erd\H os-Wintner theorems}
\bigskip\medskip
\centerline{GŽrald Tenenbaum \& Johann Verwee}
\bigskip\bigskip
{\eightpoint\leftskip1cm\rightskip1cm
\noi{\bf Abstract.} The classical theorem of Erd\H os \& Wintner furnishes a criterion for the existence of a limiting distribution for a real, additive arithmetical function. This work is devoted to providing an effective estimate for the remainder term under the assumption that the conditions in the criterion are fulfilled. We also investigate the case of a conditional distribution.
\PAR
\medskip\noi
{ \bf Keywords:}  distribution of real additive functions, mean values of complex multiplicative function, Erd\H os-Wintner theorem, effective averages, number of prime factors.\par
\smallskip 
\noi \bf 2010 Mathematics Subject Classification: \rm   11N25, 11N37, 11N60.\par }
\par 
\bigskip\bigskip
\paraun{Introduction and statement of results}
The classical theorem of Erd\H os \& Wintner \citer{E35-7-8}, \citer{EW39}, is the analogue in probabilistic number theory of Kolmogorov's  three series theorem in probability theory. It asserts that a real, additive arithmetical function $f$ possesses a limiting distribution if, and only if, the following series converge
$$\sum_{p\in\P}{\min\big(1,f(p)^2\big)\over p},\quad\sum_{\di{p\in\P}{|f(p)|\leqslant 1}}{f(p)\over p},\eqdef{3S}$$
where, here and in the sequel, $\P$ denotes the set of primes. Moreover it follows from a theorem of LŽvy \citer{PL31} that the limit law is continuous if, and only if,
$$\sum_{f(p)\neq0}{1\over p}=\infty,\eqdef{CC}$$
while a well-known theorem of Jessen and Wintner \citer{JW35} tells us that this limit law is necessarily pure. See, \eg, \citeplus{Te15}{ch. III.4} for proofs and historical comments. \par 
In this work, our first aim is to exploit a recent result of the first author \citer{Te17} on mean values of complex multiplicative functions in order to provide an effective version of the Erd\H os--Wintner theorem, or, in other words, to furnish an effective estimate for the supremum norm
$$\no{F_x-F}_\infty:=\sup_{y\in\r}\abs{F_x(y)-F(y)}\qquad (x\geqslant 1)$$
where, for each $x\geqslant 1,$ $$F_x(y):={1\over x}\sum_{\di{n\leqslant x}{f(n)\leqslant y}}1\qquad (y\in\r)$$
is the empirical distribution function and $F$ is the limiting distribution. It is well known that $F$ has characteristic function
$$\varphi_F(\tau):=\int_\r\e^{i\tau y}\d F(y)=\prod_{p}\Big(1-{1\over p}\Big)\sum_{\nu\geqslant 0}{\e^{i\tau f(\pnu)}\over \pnu}\qquad (\tau\in\r).\eqdef{fc}$$
\par 
We state our results in this direction as two separate theorems, corresponding respectively to the discrete and the continuous case. \par 
Let us first consider the situation when \eqref{3S} is realised but \eqref{CC} is not. We then define a multiplicative function $u_f$  by 
its values on prime powers $$u_f(\pnu):=\cases{1& if $f(\pnu)\neq0$,\cr
0 & if $f(\pnu)=0$,\cr}\eqdef{defuf}$$
and, given a prime $p\in\P$, write
$$S_p=S_p(f):=\sum_{\nu\geqslant 1}{u_f(\pnu)\over \pnu},\quad w_p=w_p(f):=\Big(1-{1\over p}\Big)S_p(f),\eqdef{defSw}$$
so that the convergence of the series on the left-hand side of \eqref{CC} implies the absolute convergence of $\sum_p{w_p}$.  We also plainly have
$$\alpha_f(y):=\sum_{p>y}{u_f(p)\over p}\to0,\quad\beta_f(y):={1\over \log y}\int_1^y{\alpha_f(t)\over t}\d t\to0\quad(y\to\infty).\eqdef{defabf}$$
Writing $$h_f(m):=u_f(m)\prod_{p|m}{1-1/p\over 1-w_p}\quad(m\geqslant 1),$$
we easily check that 
$$F(y):=\prod_{p}(1-w_p)\sum_{f(m)\leqslant y}{h_f(m)\over m}\qquad (y\in\r)$$
is a distribution function, indeed
$$\sum_{m\geqslant 1}{h_f(m)\over m}=\prod_{p}\Big(1+{1-1/p\over 1-w_p}S_p\Big)=\prod_{p}\Big(1+{w_p\over 1-w_p}\Big)=\prod_{p}{1\over 1-w_p}\cdot$$
With these notations, we can state our first result. Here and in the sequel, we let $\log_k$ denote the $k$-fold iterated logarithm.
\Propt{cd}{Let $f$ be a real additive function satisfying \eqref{3S} but not \eqref{CC}. Then, uniformly for $x\geqslant 2$, we have
$$\no{F_x-F}_\infty\ll R_x:=\alpha_f\big(x^{1/\log_2x}\big)+\beta_f\big(\sqrt{x}\big)^{1/4}+{1\over (\log x)^{1/6}}\cdot$$}
\noi{\it Examples.} 
(i) Let $\kappa>0$ be a parameter and consider an additive function $f$ such that $f(p)=1$ if $2^n<p\leqslant 2^n(1+1/(\log n)^\kappa)$ for some $n\geqslant 3$ and $f(p)=0$ otherwise. In this setting, the limit law is atomicÑ--\ie \eqref{CC} fails---if, and only if, $\kappa>1$. We then have
$$\alpha_f(y)\asymp \beta_f(y)\asymp {1\over (\log_2y)^{\kappa-1}},\quad R_x\asymp {1\over (\log_2x)^{(\kappa-1)/4}}\cdot$$
\par 
(ii) Assume now that $f(p)=1$ if $2^n<p\leqslant 2^n(1+1/n^\kappa)$ for some $n\geqslant 1$, while $f(p)=0$ otherwise. Then the series  \eqref{CC} converges for all $\kappa>0$ and we have
$$\alpha_f(y)\asymp{1\over (\log y)^\kappa},\quad\beta_f(y)\asymp{(\log_2y)^{\delta_{1\kappa}}\over (\log y)^{\min(1,\kappa)}},\quad R_x\asymp{1\over (\log x)^{\min(2/3,\kappa)/4}}, $$
with Kronecker's notation $\delta_{1\kappa}.$
\par 
\goodbreak
(iii) 
More generally, when the non-zero values of $f(p)$ are distributed with sufficient regularity, a simple criterion for the continuity of the limit law may be stated and subsequent estimates for $R_x$ may then be easily computed.  Indeed,  writing $$\{p\in\P:f(p)\neq0\}=\P\cap\big(\cup_{k\geqslant 1}]a_k,b_k]\big)$$ where the $a_k$, $b_k$ are integers, $2\leqslant a_k<b_k$, we first observe that this set is certainly infinite provided $$b_k>a_k+a_k^{1-c}\qquad (k\geqslant 1)\eqdef{Hoheisel}$$ for sufficiently small, positive $c$: it follows from \citer{BHP01} that, with $c=0.475$, we have $\pi(x+y)-\pi(x)\asymp y/\log x$ for $x^{1-c}\leqslant y\leqslant x$---the sharpest estimate of Hoheisel type to date. Appealing to this result and to the \pnt\ in the form
$$\sum_{a<p\leqslant b}{1\over p}=\log\Big({\log b\over \log a}\Big)+O\Big(\e^{-\sqrt{\log a}}\Big)\qquad (b\geqslant a\geqslant 3),$$
it is a simple matter to deduce that, assuming \eqref{Hoheisel}, condition \eqref{CC}  holds if, and only if, $$\sum_{k\geqslant 1}\log \Big({\log b_k\over \log a_k}\Big)=\infty.$$
\medskip
We next turn our attention to the case when \eqref{3S} and \eqref{CC} are both satisfied, which implies that the limiting distribution $F$ is continuous. We then let $\eta_f(y)$ denote any continuous, non-increasing function  tending to 0 at infinity and such that
$$\Vbs{\sum_{\di{p>y}{|f(p)|\leqslant 1}}{f(p)\over p}}\leqslant  \eta_f(y),\quad \sum_{\pnu>y}{\min(1,f(\pnu)^2)\over \pnu}\leqslant  \eta_f(y)\qquad (y\geqslant 1).\eqdef{defetay}$$
\par 
For $x\geqslant 2$, we consider a quantity $\varepsilon_x$ such that $1/\sqrt{\log x}< \varepsilon_x=o(1)$, and assume henceforth that $\varepsilon_x$ approaches 0 so slowly that
$$\eta_f\big(x^{\varepsilon_x}\big)=o\Big(\varepsilon_x^{1/3}\Big)\qquad (x\to\infty).\eqdef{Tx}$$
\par 
We write furthermore $$B_f(v)^2:=2+\sum_{\pnu\leqslant v}{f(\pnu)^2\over \pnu}\qquad (v\geqslant 1),$$
and let $\ell\mapsto Q_F(\ell):=\sup_{y\in\r}\big\{F(y+\ell)-F(y)\big\}$ denote the concentration function associated to $F$.
Since $F$ is continuous, we know that $Q_F(\ell)\to0$ as $\ell\to 0$. Effective upper bounds, depending explicitly on the sequence $\{f(p)\}_{p\in\P}$ or on $\varphi_F$ are available in the literature: see, \eg, \citer{BT12}, \citer{EK79}, \citer{Ko14}, \citer{Ru80}, and \citeplus{Te15}{ch. III.2}. For instance, the Kolmogorov--Rogozin inequality implies
$$Q_F(\ell)\ll{1\over\sqrt{1+\sum_{|f(p)|>\ell}1/p}},\eqdef{majQ1}$$
while a simple computation (see, \eg, \citeplus{Te15}{lemma III.2.9}) provides
$$Q_F(\ell)\ll\ell\int_{-1/\ell}^{1/\ell}|\varphi_F(\tau)|\d \tau.\eqdef{majQ2}$$
\Propt{cc}{Uniformly for all  real additive functions $f$ satisfying \eqref{3S} and \eqref{CC}, and all $R\geqslant 3,\,T\geqslant 3,\,x\geqslant 3$, such that $$ 2\log_2R+\dm T^2\eta_f(R)+7\leqslant \ft14\log (1/\varepsilon_x),\quad T^2\eta_f(x^{\varepsilon_x})\ll \varepsilon_x^{1/3},\eqdef{condRT}$$ we have
$$\no{F_x-F}_\infty\ll Q_F\Big({1\over T}\Big)+\varepsilon_x^{1/6}\log\Big({TB_f(R)\over \varepsilon_x}\Big)+\eta_f(R).\eqdef{majFxF}$$ }
\smallskip\goodbreak
\rems (i)  When, for instance, $\xi>1$ and $f$ is the strongly additive function defined by \hbox{$f(p)=1/(\log p)^\xi$}, an estimate of Koukoulopoulos \citer{Ko14} sharpening a result of La Bretche \& Tenenbaum \citer{BT12} yields $Q_F(\ell)\asymp \ell^{1/\xi}$ $(0<\ell\leqslant 1/3)$. Then, $B_f(R)\asymp 1$, $\eta_f(y)\asymp1/(\log y)^\xi$, the choice
$$\varepsilon_x=2/\sqrt{\log x},\quad R=\e^{c(\log x)^{1/16}},\quad T=(\log x)^{\xi/32}$$  is admissible for suitably small $c>0$, and we get, ignoring some negative powers of $\log_2 x$, $$\no{F_x-F}_\infty\ll {1\over (\log x)^{1/32}}\cdot$$
(ii) The general estimate \eqref{majFxF} may be superseded by specific known results when $f(p)$ shows rapid and smooth decrease. For instance, if $f(p)=1/p^\xi$ with $\xi>0$, $f(\pnu)=0$ $(\nu\geqslant 2)$, we have $Q(\ell)\asymp 1/|\log \ell|$ $(0<\ell\leqslant 1/3)$ by \citeplus{BT12}{Cor. 1.3}. The optimal choice is then  $$\varepsilon_x\asymp 1/\sqrt{\log x}, \quad R=\e^{c(\log x)^{1/16}},\quad\log T\asymp (\log x)^{1/24},$$ and we only get
$$\no{F_x-F}_\infty\ll{1\over (\log x)^{1/24}},$$
while the left-hand side is actually $\ll(\log_2 x)/\{(\log x)\log_3x\}$, in view of \citeplus{BT12}{Cor. 1.5}, which includes the cases of the additive functions $\log \{n/\varphi(n)\}$ and $\log \{\sigma(n)/n\}$ where $\varphi$ is Euler's totient and $\sigma(n):=\sum_{d|n}d$. This lack of precision may be traced back to the use of the general upper bounds \eqref{majtriv} and \eqref{maj2} {\it infra}, which only integrate partial information on the distribution of the $f(p)$: when $f(p)$ is quickly decreasing, a direct bound for the difference of the characteristic functions  furnishes the stated sharpening.
\par \medskip
The technique involved in the proofs of the above results is actually fairly flexible. As an illustration, we present a further effective theorem, describing how the distribution of an additive function fluctuates when restricting the support to integers with a fixed number of prime factors. To avoid technicalities we focus on the case of a strongly additive function with continuous distribution, but a completely general statement could be achieved by the same method. \par 
Let $\omega(n)$ denote the number of distinct prime factors of an integer $n$ and, for $x\geqslant 1$,  let $\pi_k(x)$ represent the cardinality of the level set $\E(x;k):=\{n\leqslant x:\omega(n)=k\}$. We have, classicaly (see \eg\ \citeplus{Te15}{ch. II.6})
$$\pi_k(x)\asymp {x(\log_2x)^{k-1}\over (k-1)!\log x}\qquad (1\leqslant k\ll\log_2x),$$
and we may replace $k-1$ by $k$ when $k\asymp\log_2x$. \par Given the strongly additive function $f$ satisfying \eqref{3S},  we consider for each $r>0$ the characteristic function
$$\varphi(\tau;r):=\prod_{p}\Big(1+{r\e^{i\tau f(p)}\over p-1}\Big)\Big(1+{r\over p-1}\Big)^{-1}\qquad (\tau\in\r),\eqdef{phitaur}$$
and denote as $\FF_r$ the corresponding distribution function. \par\goodbreak
Our estimate depends on the function $\eta_f$ defined in \eqref{defetay}. We furthermore introduce parameters $v$, $T$ and $R$ such that
 $$\eqalign{&{1\over \log_2 x}\leqslant  v\leqslant c_0,\quad  3\leqslant R\leqslant \e^{1/v},\quad T\geqslant 1,\cr
 & T^2\eta_f(R)\leqslant \log (1/v),\quad T^{2}\eta_f(x^w)\ll w\quad(w:=v^{c_1}),\cr}\eqdef{condplus}$$
where $c_0$ and $c_1$ denote strictly positive constants, depending at most on $\kappa$, $c_0$ being sufficiently small and $c_1$ sufficiently large.  \goodbreak
\Propt{T3}{Let $\kappa\in]0,1[$ and let $f$ be a real, strongly additive  function. Assume~\eqref{3S} and \eqref{CC} hold. Then, uniformly for $\kappa\leqslant r:=k/\log_2x\leqslant1/\kappa$, $y\in\r$, and $v,\,T,\,R$ satisfying~\eqref{condplus}, we have
$${1\over \pi_k(x)}\sum_{\di{n\in\E(x;k)}{f(n)\leqslant y}}1=\FF_r(y)+O\big(\gR\big)\eqdef{repEk}$$
with $$\gR:=Q_{\FF_r}\Big({1\over T}\Big)+\Big(v+{\log (1/v)\over \sqrt{k}}\Big)\log \Big({TB_f(R)\over v}\Big)+\eta_f(R)^{r/(r+1)}.$$}
\medskip
 Due to the generality of the hypotheses, this statement turns out as rather technical. Indeed an optimal choice of the parameters heavily depends on the sequence $\{f(p)\}_{p\in\P}$. However, an explicit estimate easily follows in non-pathological situations.   As an example, consider the case when $f(p):=1/(\log p)^\xi$ with $0<\xi<r$. It is then easy to show (see, \eg, \citeplus{TW14}{Exercise 259}) that $|\varphi(\tau;r)|=|\tau|^{-r/\xi}(\log |\tau|)^{O(1)}$ as $|\tau|\to\infty$ and hence, by \eqref{majQ2}, that $Q_{\FF_r}(\ell)\ll \ell$ as $\ell\to0$. We may therefore select
$$v:=1/\log_2 x,\quad  R:=\log x, \quad T:=(\log_2 x)^{\xi/2}, $$
and infer that
$\gR\ll(\log_3x)^2/\sqrt{\log_2x}+1/(\log_2x)^{\xi \min\{1/2,r/(r+1)\}}.$
\medskip
\bigskip
\paraun{The key argument}
Our approach rests on the following recent result of the first author \citeplus{Te17}{th.~1.2}, for the statement of which we introduce further notation.
We let $\M(A,B)$ designate the class of those complex-valued multiplicative functions $g$ such that
$$\max_{p}|g(p)|\leqslant A,\quad  \sum_{p,\,\nu\geqslant 2}{|g(p^\nu)|\log p^\nu\over p^\nu}\leqslant B,\eqdef{C0}
$$
and, for $\gb\in\r$, we write
$$\beta=\beta(\gb,A):=1-{\sin(2\pi\gb/A)\over 2\pi\gb/A}\cdot\eqdef{delta0}$$
Moreover, given any complex-valued function $g$, we put $\gc_g:=1$ if $g$ is real, $\gc_g:=2$ otherwise, and consider
$$M(x;g):=\sum_{n\leqslant x}g(n),\qquad Z(x;g):=\sum_{p\leqslant x}{g(p)\over p}\cdot$$
\par 
\propt{th14}{ (\citer{Te17})}{Let $$\eqalign{&\ga\in]0,\ft14], \quad\gb\in[\ga,\ft12[, \quad \gh:=(1-\gb)/\gb, \quad A\geqslant 2\gb,\quad B>0, \quad \beta:=\beta(\gb,A),\cr & 2\gb\leqslant \varrho\leqslant A,\quad x\geqslant 2, \quad 1/\sqrt{\log x}<\varepsilon\leqslant \dm,\cr}$$ and let the multiplicative functions $g$, $r$, such that \hbox{$r\in\M(2A,B)$}, $|g|\leqslant r$,  satisfy the conditions $$\leqalignno{\sum_{p\leqslant x}{r(p)-\re g(p)\over p}&\leqslant \ft12\beta\gb\log (1/\epsilon), &\eqdef{pqr}\cr
\sum_{x^\varepsilon<p\leqslant y}{\{r(p)-\re g(p)\}^\gh\log p\over p}&\ll \varepsilon^{\gc_g\delta\gh}\log y\qquad (x^\varepsilon<y\leqslant x), &\eqdef{cdpr1}\cr
\sum_{p\leqslant y}{(r(p)-\varrho)\log p\over p}&\ll \varepsilon\log y\qquad (x^\varepsilon< y\leqslant x)&\eqdef{condrho}
}$$ 
 with
$\delta\in]0,2\beta\gb/(3\gc_g)\big]$.\par \goodbreak
We then have
$$M(x;g)={\e^{-\gamma\varrho}x\over \Gamma(\varrho)\log x}\bigg\{\prod_{p}\sum_{\pnu\leqslant x}{g(\pnu)\over \pnu}+O\Big(\varepsilon^{\delta}\e^{Z(x;g)}\Big)\bigg\},\eqdef{Mf/Mr}$$
where $\gamma$ denotes Euler's constant. The implicit constant in \eqref{Mf/Mr} depends at most upon $A$, $B$, $\ga$, and $\gb$.}
\medskip
\paraun{Proof of \ref{cd}}
Let $u_f$ be defined by \eqref{defuf} and let $v_f$ be the multiplicative function defined of prime powers by $ v_f(\pnu):=1-u_f(\pnu)$. Then any integer $n\geqslant 1$ may be uniquely represented as a product $n=md$ with $u_f(m)=v_f(d)=1$, $(m,d)=1$ and $f(n)=f(m)$. Therefore
$$F_x(y)={1\over x}\sum_{\di{m\leqslant x}{f(m)\leqslant y}}u_f(m)V_m\Big({x\over m}\Big)\eqdef{dec}$$
with
$$V_m(t):=\sum_{\di{d\leqslant t}{(d,m)=1}}v_f(d)=\sum_{d\leqslant t}v_f(d;m)\qquad (t\geqslant 1), $$
say, where, for each $m$, $v_f(d;m)$ is the multiplicative function of $d$ defined on prime powers by $v_f(\pnu;m)=v_f(\pnu)$ if $p\nmid m$ and $=0$ otherwise. For $t\geqslant 3$, $m\leqslant t$,  we have, by the convergence of \eqref{CC},
$$\sum_{p\leqslant t}{1-v_f(p;m)\over p}\leqslant \sum_{p|m}{1\over p}+O(1)\leqslant \log_3t+O(1)$$
and, similarly,
$$\eqalign{\sum_{p\leqslant t}{\{1-v_f(p;m)\}\log p\over p}&\leqslant \sum_{p\leqslant t}{\{1-v_f(p)\}\log p\over p}+\sum_{p|m}{\log p\over p}\cr&=\int_1^t\sum_{u<p\leqslant t}{u_f(p)\over p}{\dd u\over u}+O(\log_2t)\leqslant \beta_f(t)\log t+O(\log_2t), \cr}$$
with  notation \eqref{defabf}.\par 
Hence, for $\varepsilon:=\beta_f(t)^{3/4}+1/\sqrt{\log t}$, we have, for $t^\varepsilon<s\leqslant t$,
$$\sum_{t^\varepsilon<p\leqslant s}{\{1-v_f(p;m)\}\log p\over p}\ll\beta_f(t)^{1/4}\log s+O(\log_2s)\ll\varepsilon^{1/3}\log t.$$
We may therefore estimate $V_m(t)$, uniformly in $m\leqslant t$, by applying \ref{th14} to $g:=v_f(\cdot;m)$ with 
$$\gb=\dm,\quad \varrho=A=\gc_g=\beta=\gh=1,\quad \varepsilon:=\beta_f(t)^{3/4}+1/\sqrt{\log t},\quad\delta=\ft13.$$  We get, for $1\leqslant m\leqslant t$,
$$V_m(t)=\Big\{1+O\Big({1\over \log 2t}\Big)\Big\}t\psi_m(t)+O\Big(t\beta_f(t)^{1/4}+{t\over (\log 2t)^{1/6}}\Big),$$
with
$$\eqalign{\psi_m(t)&:=\prod_{p\leqslant t}\Big(1-{1\over p}\Big)\prod_{\di{p\leqslant t}{p\,\nmid\, m}}\Big({1\over 1-1/p}-S_p+O\Big({1\over t}\Big)\Big)\cr
&=\Big\{1+O\Big({1\over \log 2t}\Big)\Big\}\prod_{p\leqslant t}\big(1-w_p\big)\prod_{p|m}\Big({1-1/p\over 1-w_p }\Big),\cr}$$
where $S_p$, $w_p$, are defined in \eqref{defSw} and we have taken into account that $1/(1-1/p)-S_p\geqslant 1$.
Since $w_p\leqslant 1/p$, we have $\log (1-w_p)\geqslant -2w_p$, whence
$$\prod_{p>t}\big(1-w_p\big)\geqslant \exp\Big\{-2\sum_{p>t}w_p\Big\}\geqslant 1-O\big(\alpha_f(t)\big),$$
where $\alpha_f$ is defined in \eqref{defabf}. This yields
$$u_f(m)V_m(t)=\prod_{p}(1-w_p) th_f(m)+O\big(tu_f(m)R_0(t)\big)\qquad (t\geqslant m\geqslant  1)$$
with $\dsp R_0(t):=\alpha_f(t)+\beta_f(t)^{1/4}+{1/(\log 2t)^{1/6}}.$
Splitting the sum in \eqref{dec} at $m=\fl{\sqrt{x}}$ and considering that $V_m(t)\leqslant t$, we readily obtain, uniformly for $y\in\r$,
$$F_x(y)=F(y)+O\big(E_1+E_2\big)$$
with $$\eqalign{&E_1:=\sum_{m>\sqrt{x}}{u_f(m)\over m},\quad E_2:=R_0\big(\sqrt{x}\big)\sum_{m\leqslant \sqrt{x}}{u_f(m)\over m}\ll R_0\big(\sqrt{x}\big),\cr}$$
where the last bound follows from the convergence of the series in \eqref{CC}.\par 
In order to bound $E_1$, we introduce a parameter $T\geqslant 2$ and split the summation according to whether the largest prime factor of $m$, say $P^+(m)$, exceeds  $T$ or not. We obtain, for any $\sigma\in]0,\ft13[$
$$\eqalign{E_1&\leqslant \sum_{\di{m>\sqrt{x}}{P^+(m)\leqslant T}}{1\over m}+\sum_{p>T}\sum_{\nu\geqslant 1}{u(\pnu)\over \pnu}\sum_{m\geqslant 1}{u_f(m)\over m}\cr
&\ll {1\over x^{\sigma/2}}\prod_{p\leqslant T}\Big(1+{1\over p^{1-\sigma}}\Big)+\alpha_f(T)+{1\over T}\cr}
$$
For large $T$, we select $\sigma:=4/\log T$. The last $p$-product is then $\ll\log T$, and so
$$E_1\ll x^{-2/\log T}\log T+\alpha_f(T)+{1\over T}.$$
The required estimate follows by selecting $T:=x^{1/\log_2x}.$
\bigskip
\medskip 
\paraun{Proof of \ref{cc}}

Given $R\geqslant 3$, we define the additive function $f_R$ by 
$$f_R(\pnu):=\cases{f(\pnu) & if $\pnu\leqslant R$ or $|f(\pnu)|\leqslant 1$,\cr
0 & in all other cases.\cr}\eqdef{deffR}$$
Denote by $F_x(y;R)$ the \fr\ of $f_R$ on the set of integers not exceeding $x$ and by $F(y;R)$ that of the limit law. We first observe that, when $x\in\N^*$, 
$$\abs{F_x(y;R)-F_x(y)}\leqslant \sum_{\di{\pnu>R}{|f(\pnu)|>1}}{1\over \pnu}\leqslant  \eta_f(R)\qquad (y\in\r),\eqdef{diffR}$$
the same bound being valid for $\abs{F(y;R)-F(y)}$ since the inequality is independent of $x$.
We may hence restrict to evaluating $F_x(y;R)-F(y;R)$ with the perspective of ultimately optimising the parameter $R$.
\par \goodbreak
Note that, for $3\leqslant R\leqslant x$, 
$$\eqalign{&\sum_{\pnu\leqslant x}{f_R(\pnu)\over \pnu}=\sum_{\pnu\leqslant R}{f(\pnu)\over \pnu}+\sum_{\di{R<\pnu\leqslant x}{|f(\pnu)|\leqslant 1}}{f(\pnu)\over \pnu}\ll B_f(R)\sqrt{\log_2R},\cr
&\sum_{\pnu\leqslant x}{f_R(\pnu)^2\over \pnu}=\sum_{\pnu\leqslant R}{f(\pnu)^2\over \pnu}+\sum_{\di{R<\pnu\leqslant x}{|f(\pnu)|\leqslant 1}}{f(\pnu)^2\over \pnu}\ll B_f(R)^2,\cr}$$
where we used \eqref{3S} to bound the last sum.
By the Tur‡n-Kubilius inequality, it follows, still for $3\leqslant R\leqslant x$, that
$$\eqalign{{1\over x}\sum_{n\leqslant x}\big(\e^{i\tau f_R(n)}-1\big)&={i\tau\over x}\sum_{n\leqslant x}f_R(n)+O\Big({\tau^2\over x}\sum_{n\leqslant x}f_R(n)^2\Big)\cr
&\ll |\tau|B_f(R)\sqrt{\log_2R}+\tau^2 B_f(R)^2\log_2R.\cr}\eqdef{est-tau-petit}$$
Writing
$$\varphi_x(\tau;R):={1\over x}\sum_{n\leqslant x}\e^{i\tau f_R(n)},\qquad \varphi(\tau;R):=\int_\r\e^{i\tau y}\d F(y;R)=\prod_{p}\Big(1-{1\over p}\Big)\sum_{\nu\geqslant 0}{\e^{i\tau f_R(\pnu)}\over \pnu},$$
and considering that the upper bound in \eqref{est-tau-petit} does not depend on $x$,
we hence see that
$${\varphi_x(\tau;R)-\varphi(\tau;R)\over \tau}\ll B_f(R)\sqrt{\log_2R}+|\tau| B_f(R)^2\log_2R \qquad (\tau\in\r).\eqdef{taupetit}$$
This estimate will be used for dealing with small values of $|\tau|$.
\par 
Next we evaluate $\varphi_x(\tau;R)$ when $|\tau|$ is not too close to 0, $|\tau|\leqslant T$, and assuming~\eqref{condRT}. We have, for large $x$,
$$\eqalign{\sum_{p\leqslant x}{1-\cos(\tau f_R(p))\over p}&\leqslant \sum_{p\leqslant R}{2\over p}+\sum_{\di{R<p\leqslant x}{|f(p)|\leqslant 1}}{\tau^2f(p)^2\over 2p}\cr&\leqslant 2\log_2R+7+\dm T^2\eta_f(R)\leqslant \ft14\log (1/\varepsilon_x),\cr}\qquad (\tau\in\r),\eqdef{majtriv}$$
(where we used the estimate $\sum_{p\leqslant y}1/p\leqslant \log_2y+7/2$ $(y\geqslant 2)$ which follows by partial summation from Mertens' first theorem in the form given for instance in \citeplus{Te15}{th.~I.1.8}) and similarly, for  $|\tau|\leqslant T$, since $z_x:=x^{\varepsilon_x}\geqslant R$ by \eqref{condRT},
$$\eqalign{\sum_{z_x<p\leqslant y}{\{1-\cos(\tau f_R(p))\}\log p\over p}&\leqslant \dm\tau^2\sum_{\di{z_x<p\leqslant y}{|f(p)|\leqslant 1}}{f(p)^2\over p}\log p\cr&\ll T^2\eta_f(z_x)\log y\ll\varepsilon_x^{1/3}\log y\cr}\qquad (x^{\varepsilon_x}<y\leqslant x).\eqdef{maj2}$$
\par 
We may hence apply \ref{th14} to $g:=\e^{i\tau f_R}$, with $A=\varrho=1$, $\gb=\ft12$, $\gh=1$, $\beta=1$,  $r=\1$, and $\varepsilon=\varepsilon_x$. This yields
$$\eqalign{\varphi_x(\tau;R)&=\prod_{p\leqslant x}\Big(1-{1\over p}\Big)\sum_{\pnu\leqslant x}{\e^{i\tau f_R(\pnu)}\over \pnu}+O\big(\varepsilon_x^{1/6}\big)\cr&=\prod_{p\leqslant x}\Big(1-{1\over p}\Big)\sum_{\nu\geqslant 0}{\e^{i\tau f_R(\pnu)}\over \pnu}+O\big(\varepsilon_x^{1/6}\big),\cr}$$
where we used the inequality $|\prod_{p}(u_p+v_p)-\prod_{p}u_p|\leqslant \sum_{p}|v_p|$, valid for all $u_p$, $v_p$ such that $|u_p|\leqslant 1$, $|u_p+v_p|\leqslant 1$. 
\par \goodbreak
Since
$$\eqalign{\prod_{p> x}\Big(1-{1\over p}\Big)\sum_{\nu\geqslant 0}{\e^{i\tau f_R(\pnu)}\over \pnu}&=\prod_{p> x}\Big(1+{\e^{i\tau f_R(p)}-1\over p}+O\Big({1\over p^2}\Big)\Big)\cr
&=\exp\Bigg\{\sum_{\di{p>x}{|f(p)|\leqslant 1}}{i\tau f(p)\over p}+O\Big({\tau^2f(p)^2\over p}\Big)+O\Big({1\over x\log x}\Big)\Bigg\}\cr
&=\exp\Big\{O\Big(\eta_f(x)(1+\tau^2)+{1\over x\log x}\Big)\Big\},\cr}$$
(using $|\tau|\leqslant 1+\tau^2$) we eventually obtain, for $|\tau|\leqslant T$,
$$\varphi_x(\tau;R)=\varphi(\tau;R)\Big\{1+O\Big(\varepsilon_x^{1/6}\Big)\Big\}+O\big(\varepsilon_x^{1/6}\big),$$
and so
$$\varphi_x(\tau;R)=\varphi(\tau;R)+O\big(\varepsilon_x^{1/6}\big).\eqdef{estphixR}$$
This enables an appeal to the Berry-Esseen inequality
$$\no{F_x(\cdot;R)-F(\cdot;R)}_\infty\ll Q\Big({1\over T};R\Big)+\int_{-T}^T\vbs{{\varphi_x(\tau;R)-\varphi(\tau;R)\over \tau}}\d\tau,$$
where $Q(\cdot;R)$ is the concentration function associated to $F(\cdot;R)$. Taking \eqref{diffR} into account, we get
$$\no{F_x-F}_\infty\ll Q_F\Big({1\over T}\Big)+\int_{-T}^{T}\vbs{{\varphi_x(\tau;R)-\varphi(\tau;R)\over \tau}}\d\tau+\eta_f(R).$$
To bound the last integral, say $I$, from above, we introduce a parameter $u\in]0,1[$ and apply \eqref{taupetit} for $|\tau|\leqslant u$, then \eqref{estphixR} for $u<|\tau|\leqslant T$. This yields
$$\eqalign{I&\ll uB_f(R)\sqrt{\log_2R}+u^2B_f(R)^2\log_2R+\varepsilon_x^{1/6}\log (T/u)\cr
&\ll \varepsilon_x^{1/6}\log\Big({TB_f(R)\over \varepsilon_x}\Big) \cr}$$
for the quasi-optimal choice $u:=\varepsilon_x^{1/6}/\{B_f(R)\sqrt{\log_2R}\}$.
\bigskip\medskip
\paraun{Proof of \ref{T3}}
Let $f$ be strongly additive, satisfying \eqref{3S} and \eqref{CC}, and for $R\geqslant 3$ let $f_R$ be defined by~\eqref{deffR}. We start with a lemma showing that, for large $R$, we have $f_R(n)=f(n)$ for most integers $n\in\E(x;k)$. We recall the notation $r:=k/\log_2x$ and put
$$\sigma_f(R):=\eta_f(R)^{r/(r+1)}+1/(\log x)^{r/(r+1)}.$$\par 
\Propl{fR-f}{Let $\kappa\in]0,1[$. Uniformly for $\kappa\leqslant r:=k/\log_2x\leqslant 1/\kappa$, $3\leqslant R\leqslant x$, we have
$$\sum_{\di{n\in\E(x;k)}{f_R(n)\neq f(n)}}1\ll \sigma_f(R)\pi_k(x).\eqdef{estfRf}$$}
\nid We may plainly assume $x$ to be large and hence that $k\geqslant 2$. Put $$\P_R:=\{p\in\P:p>R,\,|f(p)|>1\}, \quad E_R(x):=\sum_{\di{p\leqslant x}{p\in\P_R}}{1\over p}\leqslant \eta_f(R).$$ The quantity to be bounded does not exceed the number of those integers $n\in\E(x;k)$ having at least one prime divisor in $\P_R$.  \par 
From the classical Hardy-Ramanujan estimate for $\pi_k(y)$ (see \eg\ \citeplus{Te15}{Ex. 264}) the left-hand side of \eqref{estfRf} is, for an absolute constant $a$,
$$\leqslant \sum_{n\in\E(x;k)}\sum_{\di{p^\nu\|n}{p\in\P_R}}1\ll\sum_{\di{p^\nu\leqslant x}{p\in\P_R}}\pi_{k-1}\Big({x\over p^\nu}\Big)\ll\sum_{\di{p^\nu\leqslant x}{p\in\P_R}}{x\{\log_2(3x/p^\nu)+a\}^{k-2}\over p^\nu(k-2)!\log (2x/p^\nu)}\cdot$$
Put $v:=\sigma_f(R)^{1/r}$. The subsum corresponding to $p^\nu\leqslant x^{1-v}$ is plainly
$$\ll{x(\log_2x)^{k-2}E_R(x)\over v(k-2)!\log x}\ll{\sigma_f(R)x(\log_2x)^{k-2}\over (k-2)!\log x}\ll\sigma_f(R)\pi_k(x).$$
The complementary subsum may be dealt with by partial summation. By the prime number theorem, it is
$$\eqalign{&\ll{1\over (k-2)!}\int_{x^{1-v}}^x{x\{\log_2(3x/t)+a\}^{k-2}\over t\log (2x/t)\log t}\d t\asymp{x\over (k-2)!\log x}\int_1^{x^v}{(\log_23u+a)^{k-2}\over u\log 2u}\d u\cr
&\ll{x\{\log_2x-\log (1/v)+a\}^{k-1}\over (k-1)!\log x}\ll\pi_k(x)\Big\{1-{\log (1/v)\over \log_2x}\Big\}^k\ll\pi_k(x)v^r=\sigma_f(R)\pi_k(x).\cr}$$
\vskip-3mm\qed
\medskip
Our next lemma consists in obtaining a uniform upper bound for $$S_R(x;\tau,z):=\sum_{n\leqslant x}z^{\omega(n)}\e^{i\tau f_R(n)}\qquad (x\geqslant 1,\,|z|=r).$$
\Propl{LmajSR}{Let $\kappa\in]0,1[$. Uniformly for $3\leqslant R\leqslant \log x$, $\kappa\leqslant r\leqslant 1/\kappa$, $z=r\e^{i\vartheta}$, $|\vartheta|\leqslant \pi$, $|\tau|\leqslant T$, we have
$$S_R(x;\tau,z)\ll x(\log x)^{r-1}\bigg\{{\e^{9rT^2\eta_f(R)}\over (\log x)^{r\vartheta^2/60}}+{1\over \sqrt{\log x}}\bigg\}\cdot\eqdef{majSR}$$}
\nid By \citeplus{Te17}{cor.~2.1}, the left-hand side of \eqref{majSR} is
$$\ll x(\log x)^{r-1}\bigg\{{1+m_f(x;\tau)\over \e^{m_f(x;\tau)}}+{1\over \sqrt{\log x}}\bigg\}\eqdef{majSR1}$$
with 
$$m_f(x;\tau)=r\min_{|t|\leqslant \log x}\sum_{p\leqslant x}{1-\cos\{\vartheta+t\log p+\tau f_R(p)\}\over p}.$$
Let $\no a$ denote the distance of the real number $a$ to the set of integers. The elementary inequality $\no{a+b}^2\geqslant \dm\no a^2-b^2$ and the standard lower bound $1-\cos a\geqslant 8\no{a/2\pi}^2$ yield
$$m_f(x;\tau)\geqslant 4 r\min_{|t|\leqslant \log x}\lambda_f(x;t)-8r\tau^2\eta_f(R)\eqdef{minmf}$$
 with $$\lambda_f(x,t):=\sum_{R<p\leqslant x}{1\over p}\NO{\vartheta+t\log p\over 2\pi}^2\cdot$$
 \goodbreak
Now by \citeplus{Te15}{lemma III.4.13}, we have, restricting the $p$-sum to $y<p\leqslant x$ with $y>R$,
$$\lambda_f(x;t)\geqslant \ft1{12}\log \Big({\log x\over \log y}\Big)+O\Big({1\over |t|\log y}+{1+|t|\over \e^{\sqrt{\log y}}}\Big)\qquad (2\leqslant y\leqslant x).\eqdef{valp}$$
If $1\leqslant |t|\leqslant  \log x$, we select $y:=\exp\big\{(\log_2x)^2\big\}$ to get
$$\lambda_f(x;t)\geqslant \ft1{12}\log_2x+O(\log_3x).$$
Let us then define $\nu:=(\log x)^{(\vartheta^2/2\pi^2)-1}.$ If $\nu\leqslant |t|\leqslant 1$, we select $y=\e^{1/\nu}$ in \eqref{valp} and obtain
$$\lambda_f(x;t)\geqslant {\vartheta^2\over 24\pi^2}\log_2x+O(1).$$
Finally, if $|t|\leqslant \nu$, we have
$$\eqalign{\lambda_f(x;t)&\geqslant\sum_{\log x<p\leqslant \e^{1/\nu}}{\vartheta^2/4\pi^2+O(\nu\log p)\over p}\geqslant {\vartheta^2(1-\vartheta^2/2\pi^2)\over 4\pi^2}\log\Big({\log x\over \log_2x}\Big)+O(1)\cr&\geqslant {\vartheta^2\over 8\pi^2}\log\Big({\log x\over \log_2x}\Big)+O(1).\cr}\eqdef{minlambda}$$
Carrying back into \eqref{minmf} and \eqref{majSR1} yields the stated estimate since $1/6\pi^2>1/60$.
\qed
\medskip
We now deduce from \ref{th14} an asymptotic formula with remainder for $S_R(x;\tau,z)$ when $z$ belongs to a neighbourhood of the real point $r$ on the circle $|z|=r$. 
\Propl{SR}{For suitable constant $c_0\in]0,1]$, arbitrary $c_1>0$, both depending at most on $\kappa$, and uniformly under the assumptions 
$$\eqalign{&z=r\e^{i\vartheta},\  {1\over \log_2 x}\leqslant v\leqslant c_0,\  |\vartheta|\leqslant \vartheta_x:=\sqrt{c_1\log(1/v)\over  \log_2 x},\  1+|\tau|\leqslant  T,\cr
&  3\leqslant R\leqslant \e^{1/v}, \quad T^2\eta_f(R)\leqslant \log (1/v),\quad T^{2}\eta_f(x^w)\ll w\quad(w:=v^{c_1}),\cr}\eqdef{condfin}$$
 we have
$$S_R(x;\tau,z)={x\e^{-\gamma r}\over \log x}\bigg\{\prod_{p\leqslant x}\Big(1+{z\e^{i\tau f_R(p)}\over p-1}\Big)+O\Big((|\vartheta|+v^{2})(\log x)^{r}\Big)\bigg\},\eqdef{evalSR}$$}
\nid We apply \ref{th14} with $r(n):=r^{\omega(n)}$, $g(n):=z^{\omega(n)}\e^{i\tau f_R(n)}$, $\gb:=\ft12\min(1,r)$, $A:=\max(1,r)$, $\varrho:=r$, $\delta:=c_2\beta\gb $, and $\varepsilon:=(|\vartheta|+v^{2})^{1/\delta}$. We select $c_2$ so small to ensure that $2\delta\gh\leqslant 1$, where $\gh=(1-\gb)/\gb$.
\par 
Since
$$\eqalign{\sum_{p\leqslant x}{r(p)-\re g(p)\over p}&=\sum_{p\leqslant R}{r-\re g(p)\over p}+r\sum_{\di{R<p\leqslant x}{|f(p)|\leqslant 1}}{1-\cos\big(\vartheta+\tau f(p)\big)\over p}+O(1)\cr
&\leqslant 2r\log_2R+r\vartheta^2\log_2x+rT^2\eta_f(R)+O(1)\cr
&\leqslant (c_1+3)r\log(1/v)+O(1),  \cr}$$
we see that condition \eqref{pqr} is satisfied for an appropriate choice of $c_0$ and $c_2$: indeed, this is clear if $v^{2}>|\vartheta|$ for then $1/v<\sqrt{2}/\varepsilon^{\delta/2}$, and, if $v^{2}\leqslant |\vartheta|$, we have $\varepsilon^\delta\leqslant 2\vartheta_x$ whence $\log (1/v)\leqslant \log_3x\ll \log 1/\vartheta_x\ll \delta\log (1/\varepsilon). $
\par 
Next,  since $v^{c_1}\leqslant \varepsilon$ provided $c_1\geqslant 2/\delta$, we have, for $x^\varepsilon<y\leqslant x$,
$$\sum_{\di{x^\varepsilon<p\leqslant y}{|f(p)|\leqslant 1}}{r(1-\cos\{\vartheta+\tau f(p)\})^\gh\log p\over p}\ll\big(|\vartheta|^{2\gh}+T^{2}\eta_f(x^\varepsilon)\big)\log y\llÊ\varepsilon^{2\delta\gh}\log y,$$
and so condition \eqref{cdpr1} is also satisfied. 
Considering the fact that \eqref{condrho} holds trivially, we obtain
$$S_R(x;\tau,z)={x\e^{-\gamma r}\over \log x}\bigg\{\prod_{p\leqslant x}\Big(1+{z\e^{i\tau f_R(p)}\over p-1}\Big)+O\Big(\varepsilon^\delta\e^{z\sum_{p\leqslant x}\e^{i\tau f_R(p)}/p}\Big)\bigg\}.$$
The required estimate hence follows from a trivial estimate for the last sum over $p$.\qed
\bigskip
We are now in a position to embark on the final part of the proof.
\par 
Define
$$L(\tau;x):=\sum_{p\leqslant x}{\e^{i\tau f_R(p)}\over p-1},\quad G_\tau(z;x):=\e^{-zL(\tau;x)}\prod_{p\leqslant x}\Big(1+{z\e^{i\tau f_R(p)}\over p-1}\Big).$$
Under  conditions \eqref{condfin} for $\tau$, we have 
$$\eqalign{\abs{L(0;x)-L(\tau;x)}&\leqslant 2\log_2R+T\eta_f(R)+\ft12T^2\eta_f(R)+O(1)\cr&\leqslant 4\log (1/v)+O(1)\ll\log k,\cr}\eqdef{L0-Ltau}$$
in particular $L(\tau;x)=\log_2x+O(\log k)$, while of course $L(0;x)=\log_2x+O(1)$.
 Moreover, $G_\tau(z;x)$ is an entire function of $z$ which is uniformly bounded on any compact subset with respect to $\tau$ and $x$, so we have for instance $$G_\tau^{(j)}(0;x)/j!\ll 1/(1+r)^j\quad(j\geqslant 0).\eqdef{majG}$$
\par 
We now apply Cauchy's integral formula to $S_R(x;\tau,z)$ for the circle $|z|=r=k/\log_2x$, under hypotheses \eqref{condplus}. This yields
 $$\sum_{\di{n\leqslant x}{\omega(n)=k}}\e^{i\tau f_R(n)}={x\e^{-\gamma r}(I_1-I_2+I_3)\over\log x}+I_4$$
 with
 $$\eqalign{I_1&:={1\over 2\pi i}\oint_{|z|=r}\e^{zL(\tau;x)}G_\tau(z;x){\dd z\over z^{k+1}},\quad I_2:={1\over 2\pi i}\int_{\di{|z|=r}{|\vartheta_x<|\vartheta|\leqslant \pi}}\e^{zL(\tau;x)}G_\tau(z;x){\dd z\over z^{k+1}},\cr
I_3&:={\e^{k}\over 2\pi i}\int_{\di{|z|=r}{|\vartheta|\leqslant \vartheta_x}}O(|\vartheta|+v^2){\dd z\over z^{k+1}},\quad
 I_4:={1\over 2\pi i}\int_{\di{|z|=r}{\vartheta_x<|\vartheta|\leqslant \pi}}S_R(x;\tau,z){\dd z\over z^{k+1}}\cdot\cr}
 $$
 \par 
The main term is provided by $I_1$,  equal to the coefficient of $z^k$ in $\e^{zL(\tau;x)}G_\tau(z;x)$. We thus have
$$I_1=\sum_{0\leqslant j\leqslant k}{L(\tau;x)^{k-j}G_\tau^{(j)}(0;x)\over (k-j)!j!}={L(\tau;x)^k\over k!}\Big\{G_\tau(r;x)+O\Big({\log k\over k}\Big)\Big\},\eqdef{evalsum}$$
by \eqref{L0-Ltau} and \eqref{majG}, after a short computation involving truncating the sum at $\big\lfloor\sqrt{k}\big\rfloor$, for instance, and noting that
$${1\over (k-j)! L(\tau;x)^{j}}={r^j\over k!}\Big\{1+O\Big({j^2+j\log k\over k}\Big)\Big\}\qquad \big(0\leqslant j\leqslant \sqrt{k}\big).$$\par 
 The integrals $I_j$ $(2\leqslant j\leqslant 4)$ are treated as error terms. We first have, by \eqref{L0-Ltau},
 $$\eqalign{I_2&\ll \int_{\vartheta_x}^\pi \e^{-r L(0;x)\cos \vartheta}\d\vartheta\ll {\sqrt{k}(\log_2x)^k\over v^{4r}k!}\int_{\vartheta_x}^\pi\e^{-2k\vartheta^2/\pi^2}\d\vartheta\cr
 &\ll{\e^{-k\vartheta_x^2/\pi^2}(\log_2x)^k\over v^{4r}k!}\ll{v^{(c_1/\pi^2-4)r}(\log_2x)^k\over k!}\ll{v(\log_2x)^k\over k!},\cr}$$
 since $c_1$ may be taken arbitrarily large in terms of $\kappa$.
Next
$$I_3\ll {\e^k\over r^k}\{\vartheta_x^2+v^{2}\vartheta_x\}\ll{x(\log_2x)^k\over k!\log x}\Big\{v+{\log 1/v\over \sqrt{k}}\Big\}\asymp\pi_k(x)\Big\{v+{\log 1/v\over \sqrt{k}}\Big\},$$
and finally, by \eqref{majSR}, still with a suitable choice of $c_1$,
$$I_4\ll{x\e^k\over r^k\log x}\Big\{{v^{c_1r/60}\over v^{9r}\sqrt{k}}+{1\over \sqrt{\log x}}\Big\}\ll v\pi_k(x).$$
\par 
Gathering our estimates, we arrive at
$$\sum_{\di{n\leqslant x}{\omega(n)=k}}\e^{i\tau f_R(n)}={x\e^{-\gamma r}L(\tau;x)^kG_\tau(r;x)\over  k!\log x}+ O\bigg(\pi_k(x)\Big(v+{\log (1/v)\over \sqrt{k}}\Big)\bigg),$$
since the error term of \eqref{evalsum} may be absorbed by the other remainders.
Applying this with $\tau=0$, we get
$$\pi_k(x)={x\e^{-\gamma r}L(0;x)^kG_0(r;x)\over  k!\log x}\Big\{1+O\Big(v+{\log (1/v)\over \sqrt{k}}\Big)\Big\},$$
 and so, recalling definition \eqref{phitaur},
$$\eqalign{\varphi_x(\tau;k)&:={1\over \pi_k(x)}\sum_{\di{n\leqslant x}{\omega(n)=k}}\e^{i\tau f_R(n)}\cr&={L(\tau;x)^kG_\tau(r;x)\over L(0;x)^kG_0(r;x)}\Big\{1+ O\Big(v+{\log (1/v)\over \sqrt{k}}\Big)\Big\}+O\Big(v+{\log (1/v)\over \sqrt{k}}\Big)\cr
&={L(\tau;x)^kG_\tau(r;x)\over L(0;x)^kG_0(r;x)}+O\Big(v+{\log (1/v)\over \sqrt{k}}\Big)=\varphi(\tau;r)+O\Big(v+{\log (1/v)\over \sqrt{k}}\Big),\cr}$$
in view of \eqref{L0-Ltau}
and since
$$
\abs{\sum_{p> x}{\e^{i\tau f_R(p)}-1\over p-1}}\leqslant T\eta_f(x)+\dm T^2\eta_f(x)+O\Big({1\over x}\Big)\ll v.$$
\par
It remains to apply the Berry-Esseen inequality, taking into account  a variant of \eqref{est-tau-petit} conditioned to $\omega(n)=k$ in order to handle the contribution of small $\tau$. Assuming~\eqref{condplus},
we get, for $0<u\leqslant T$,
$$\eqalign{\sup_{y\in\r}\Vbs{{1\over \pi_k(x)}&\sum_{\di{n\in\E(x;k)}{f_R(n)\leqslant y}}1-\FF_r(y)}\ll {Q_{\FF_r}\Big({1\over T}\Big)}+\int_{-T}^T\abs{\varphi_x(\tau;k)-\varphi(\tau;r)\over \tau}\d\tau\cr
&\ll{Q_{\FF_r}\Big({1\over T}\Big)}+ uB_f(R)\sqrt{\log_2R}+u^2B_f(R)^2+\Big(v+{\log (1/v)\over \sqrt{k}}\Big)\log {T\over u}\cr
&\ll{Q_{\FF_r}\Big({1\over T}\Big)}+ \Big(v+{\log (1/v)\over \sqrt{k}}\Big)\log \Big({TB_f(R)\over v}\Big),\cr}$$
with the quasi-optimal choice $u=v/\{B_f(R)\sqrt{\log_2R}\}$. 
The required estimate now follows from \eqref{estfRf}.
\bigskip
\noi{\it Acknowledgements.} The authors wish to thank the referee for a thorough checking of this work and useful suggestions. The first author takes pleasure in expressing warm thanks to RŽgis de la Bretche for useful remarks  on a preliminary draft.  \bigskip
\goodbreak
\centerline{\twelvebf References}
\bigskip
\eightpoint{
\bibtem{BHP01} R.C. Baker, G. Harman \& J. Pintz, 
The difference between consecutive primes, II,
{\it Proc. London Math. Soc.} (3) {\bf83} \numero3 (2001), 532Ð562. 
\bibtem{BT12} R. de la Bretche \& G. Tenenbaum, Sur la concentration de certaines fonctions additives, {\it Math. Proc. Camb. Phil. Soc. \bf 152}, \numero1 (2012), 179--189; erratum {\it ibid.}, 191.\par 
\bibtem{E35-7-8} P. Erd\H os, On the density of some sequences of numbers I, {\it J.
London Math. Soc.
\bf 10} (1935), 120--125; II, {\it ibid.} {\bf 12} (1937), 7--11 ; III, {\it ibid.} {\bf 13} (1938), 119--127.
\par 
\bibtem{EK79} P. Erd\H os \& I. K‡tai, On the concentration of distribution of additive functions, {\it Acta Sci. Math. (Szeged) \bf41} \numeros 3-4 (1979), 295--305. \par
\bibtem{EW39} P. Erd\H os \& A. Wintner, Additive arithmetical functions and statistical independence, {\it
Amer. J. Math. \bf 61} (1939), 713--721. \par 
\bibtem{JW35} B. Jessen \& A. Wintner, Distribution functions and the Riemann Zeta function, {\it Trans. Amer. Math.
Soc. \bf 38} (1935), 48--88.\par 
\bibtem{Ko14} D. Koukoulopoulos, 
On the concentration of certain additive functions,
{\it Acta Arith. \bf162}, no. 3 (2014),  223Ð241.\par 
\bibtem{PL31} P. LŽvy, Sur les sŽries dont les termes sont des variables Žventuelles
ind\'epen\-dantes, {\it Studia Math. \bf 3} (1931), 119--155.\par 
\bibtem{Ru80} I.Z. Ruzsa,
On the concentration of additive functions,
{\it Acta Math. Acad. Sci. Hungar.} {\bf 36} (1980), no. 3-4, 215--232.\par 
\bibtem{Te15} G. Tenenbaum, {\it Introduction to analytic and probabilistic number theory}, 3rd ed., Graduate Studies in Mathematics 163, Amer. Math. Soc. 2015.\par
\bibtem{Te17} G. Tenenbaum, Moyennes effectives de fonctions multiplicatives complexes, {\it Ramanujan J. \bf44}, no. 3 (2017), 641--701; corrigendum, {\it ibid. \bf 52} no. 4 (2020).\par
\bibtem{TW14} G. Tenenbaum \& J. Wu (coll.), {\it ThŽorie analytique et probabiliste des nombres, 307 exercices corrigŽs},   coll. ƒchelles, Belin, 2014, 347 pp.\par }
\bigskip
{\leftskip2mm\rightskip-2cm\sevenrm
\gutter=4cm \doublecolumns
 \obeylines \baselineskip=7pt
 G\'erald Tenenbaum\par
Institut \'Elie Cartan\par 
Universit\'e de Lorraine\par
 BP 70239\par
54506 Vand\oe uvre-ls-Nancy Cedex\par
 France
\smallskip
{\seventt gerald.tenenbaum@univ-lorraine.fr}
Johann Verwee
Institut ƒlie Cartan
 UniversitŽ de Lorraine
 BP 70239
 54506 Vand{\oe}uvre-ls-Nancy Cedex
 France
 \smallskip
 {\seventt johann.verwee@univ-lorraine.fr}
\singlecolumn
\par
\par}
\end